\documentclass[10pt,a4paper]{article}

\usepackage[left=20mm, right=15mm, top=15mm, bottom=15mm]{geometry}

\usepackage{subfig}
\usepackage{flushend}
\usepackage{indentfirst}
\usepackage{graphics}
\usepackage{amsmath}
\usepackage{graphicx}
\usepackage{epstopdf}
\usepackage{float}

\usepackage[font=footnotesize,labelfont=bf]{caption}
\usepackage[labelsep=period]{caption}
\usepackage{graphicx}%
\usepackage{multirow}%
\usepackage{amsmath,amssymb,amsfonts}%
\usepackage{amsthm}%
\usepackage{mathrsfs}%
\usepackage[title]{appendix}%
\usepackage{xcolor}%
\usepackage{textcomp}%
\usepackage{manyfoot}%
\usepackage{booktabs}%
\usepackage{algorithm}%
\usepackage{algorithmicx}%
\usepackage{algpseudocode}%
\usepackage{listings}%

\newtheorem{theorem}{Theorem}
\theoremstyle{thmstyletwo}%
\newtheorem{Example}{Example}%
\newtheorem{lemma}{Lemma}%
\theoremstyle{thmstylethree}%
\newtheorem{definition}{Definition}%
\graphicspath{{./figure/}}

\begin{document}

\title{\huge \textbf{Viscosity Iterative algorithm for solving Variational Inclusion and Fixed point problems involving Multivalued Quasi-Nonexpansive and Demicontractive Operators in real Hilbert Space.}}

\date{}

\begin{@twocolumnfalse}
\maketitle

\author{\textbf{Furmose Mendy}$^{1,*}$, \textbf{John T Mendy}$^{2}$\\\\
\footnotesize $^{1}${College of Mathematics and Statistics, University of Toledo, Toledo, 43607, Ohio, United States\\Email: fmendy@rockets.utoledo.edu}\\
\footnotesize $^{2}${College of Mathematics and Statistics, University of Toledo, Toledo, 43607, Ohio, United States\\Email: jmendy@rockets.utoledo.edu}\\
\footnotesize $^{*}$Corresponding Author: jmendy@rockets.utoledo.edu}\\\\\\

\end{@twocolumnfalse}

\noindent \textbf{\large{Abstract}} \hspace{2pt} {This paper presents a modified general viscosity iterative process designed to solve variational inclusion and fixed point problems involving multi-valued quasi-nonexpansive and demi-contractive operators.  The modified iterative process incorporates a viscosity approximation technique to handle the nonexpansive and contractive mappings, providing a more robust and efficient solution approach. By introducing an additional sequence of iterates, the algorithm iteratively approximates the desired solution by combining fixed point iteration with viscosity approximation. The proposed method has been proven to converge strongly to the solution of the given problem, ensuring the reliability and accuracy of the results.\\

\noindent \textbf{\large Keywords}: Iterative algorithm, Quasi-nonexpansive, Multi-valued demi-closed mapping, Hilbert space, Variational inequality, Strongly Monotone Mappings \\
\noindent \hrulefill \\
\noindent \textbf{\large Subject Classifications [2010]}: Primary 47H10; Secondary 47H09, 47J05, 47J25, 54H25.

\section{Introduction}
For several years, the study of fixed point theory for multi-valued mapping has garnered significant attention and continues to be a subject of great interest in the field of mathematics. Fixed point theory deals with the existence and properties of points within a mathematical function that do not move when the function is applied.

Multi-valued mappings, also known as set-valued mappings, are functions that assign multiple outputs to a single input. These mappings play a crucial role in various areas of mathematics, such as optimization, game theory, economics, and mathematical physics.

The study of fixed point theory for multi-valued mappings has proved to be fruitful due to its wide range of applications and its profound theoretical implications. Researchers have been exploring different aspects of this theory, including the existence and uniqueness of fixed points, stability, approximation methods, and convergence properties.

One of the main reasons for the persistence of interest in this area of study is its ability to provide solutions to a wide range of problems across different disciplines. For example, in economics, fixed point theory can help analyze equilibria in markets with multiple players, while in optimization, it can aid in finding optimal solutions for complex systems.

Moreover, fixed point theory for multi-valued mappings has also contributed to the development of related areas of research, such as operator theory, topological degree theory, and convex analysis. By investigating the properties and behavior of fixed points in multi-valued mappings, mathematicians have gained a deeper understanding of these fields and have been able to establish connections and develop new techniques.(See, {\color{blue}\cite{4},\cite{18,19},\cite{9},\cite{26},\cite{8}} and {\color{blue}\cite{23}})
\quad\\
\quad\\
Viscosity iterative algorithms have been extensively studied in recent years for finding common fixed points of single-valued nonexpansive mappings and solving variational inequality problems. These investigations have built upon the concepts of viscosity solutions introduced by various researchers. (see e.g {\color{blue}\cite{15}, \cite{JT0}, \cite{F}, \cite{14},  \cite{28},  \cite{JT2},  \cite{JT3},  \cite{JT8},  \cite{sow}})
\quad\\
\quad\\
Throughout this paper, we denote $H$ to be real Hilbert space with the inner product $\displaystyle \langle.,.\rangle$ induced by the norm $\displaystyle \|.\|$. Let $\mathbb{K}$, to be a nonempty, closed and convex subset of $H$.
\quad\\
An operator $\Lambda: H \to H$ is said to be Lipschitz if there exists a constant $L>0$ such that
\begin{equation}\label{I}
    \|\Lambda \psi - \Lambda \pi\|\leq L\|\psi - \pi\|, \forall \psi, \pi\in H
\end{equation}
$\displaystyle \Lambda: H \to H$ is said to be strongly positive if there exists a constant $k>0$ such that
\begin{equation}\label{I0}
   \langle \Lambda \psi, \psi\rangle\geq k\|\psi\|^{2}, \quad\forall \psi\in H
\end{equation}
$\displaystyle \Lambda: H \to H$ is said to be $k-$strongly monotone if there exists a constant $k\in (0,1)$ such that
\begin{equation}\label{I1}
   \langle \Lambda \psi - \Lambda \pi, \psi - \pi\rangle_{H}\geq k\|\psi - \pi\|^{2}, \quad\forall \psi, \pi\in H
\end{equation}
\begin{definition}
A multivalued mapping
\begin{enumerate}
  \item $\displaystyle T: D(T)\subseteq H\to CB(D)$ is called $L-$Lipschitzian if there exists $L>0$, such that
  \begin{equation*}
    \mathcal{H}(T\psi,T\pi) \leq L\|\psi - \pi\|, \forall \psi,\pi \in D(T)
  \end{equation*}
  and $T$ is contraction if $L\in(0,1)$ and noneaxpansive if $L = 1$.
  \item $T$ is called quasi-nonexpansive if $\displaystyle \mathcal{H}(T\psi, Tp)\leq \|\psi - p\|, \quad \forall \psi\in D(T), p\in Fix(T)$
  \item $\displaystyle T: D(T)\subseteq H \to CB(D)$ is said to be $k-$stritly pseudo-contractive, if there exists $k\in(0,1)$ such that for all $\psi,\pi\in D(T)$, the following holds;
      \begin{equation*}
        \Big(\mathcal{H}(T\psi, T\pi)\Big)^{2} \leq \|\psi - \pi\|^{2} + k\|(\psi - u) - (\pi - v)\|^{2}, \forall u\in T\psi, v\in T\pi
      \end{equation*}
  If $k =1$, the map $T$ is said to be pseudocontractive.
\item \cite{23} $\displaystyle T: D(T) \subseteq E\to 2^{E}$ is said to be demicontractive if $Fix(T) \neq \emptyset$ and for all $p\in Fix(T), \psi\in D(T)$ there exists $k \in (0,1)$ such that
    \begin{equation*}
        \Big(\mathcal{H}(T\psi, Tp)\Big)^{2} \leq \|\psi - p\|^{2} + kd(\psi, T\psi)^{2}.
      \end{equation*}
 If $k =1$, the map $T$ is said to be hemicontractive.
\end{enumerate}
\end{definition}
Let $(X, d)$ be a metric space, $\mathbb{K}$ be a nonempty subset of $X$ and $\displaystyle T : \mathbb{K} \to 2^{\mathbb{K}}$ be a multivalued
mapping. An element $\psi \in \mathbb{K}$ is called a fixed point of $T$ if $\psi \in T \psi$. The fixed point set of $T$
is denoted by $\displaystyle Fix(T) := \{\psi \in D(T) : \psi \in T \psi\}$ where $\displaystyle D(T) := \{\psi \in X : T \psi \neq \emptyset\}$. It is easy
to see that single-valued mapping is a particular case of multivalued mapping.
\quad\\
Let $D$ be a nonempty suset of a normed  linear space $E$. The set $D$ is called proximinal (see\cite{20}) if for each $\psi\in E$, there exists $u\in D$ such that
\begin{equation}\label{X3}
    d(\psi,u): = \inf\{\|\psi - \pi\|: \pi\in D\}, \forall \psi,\pi\in E
\end{equation}
where $\displaystyle d(\psi,\pi): = \|\psi - \pi\|$ for all $\psi,\pi\in E$. Every closed, nonempty and convex set of real Hilbert space is proximinal. The family of nonempty closed bounded subsets, nonempty compact subsets, and nonempty proximinal bounded subsets be donated as $CB(D), K(D)$ and $P(D)$ respectively.
\quad\\
Let $\Lambda, B \in CB(D)$. Then the Hausdorff metric in $\mathcal{H}$ id defined by
\begin{equation}\label{X4}
    \mathcal{H}(\Lambda,B) = \max\Big\{\sup_{a\in \Lambda} d(a,B), \sup_{b\in B}d(b,\Lambda)\Big\}.
\end{equation}

Let $\displaystyle \Lambda : D(\Lambda) \subset H \to 2^{H}$ be a multivalued operators. Then, recall that  $\Lambda$ is monotone if $\langle u - v, \psi - \pi\rangle \ge0, (\psi,u),(\pi,v)\in G(\Lambda)$ such that
\begin{equation}\label{X}
    G(\Lambda): = \{(\psi,u): \psi \in D(\Lambda), u \in \Lambda \psi\}
\end{equation}
A monotone mapping $\displaystyle \Lambda : H \to 2^{H}$ is said to be maximal if its graph $G(A)$ is not properly
contained in the graph of any other monotone mapping.
\quad\\
A mapping $\displaystyle \Lambda: H\to H$ is said to be $\alpha-$inverse strongly monotone if there exits a constant $\alpha >0$ such that
\begin{equation}\label{X1}
    \langle \Lambda \psi - \Lambda \pi, \psi - \pi\rangle_{H}\geq \alpha\|\Lambda \psi - \Lambda \pi\|^{2}, \quad \forall \psi,\pi\in H
\end{equation}
%

Let $\displaystyle \Lambda: H \to H$ be a single-valued nonlinear mapping and $\displaystyle \Pi: H \to 2^{H}$ be a set-valued mapping. Then the variational inclusion problem is as follows: Find $\psi\in H$, such that
\begin{equation}\label{X2}
    \omega \in \Pi(\psi) + \Lambda(\psi)
\end{equation}
where $\omega$ is the zero vector in $H$. We denote the solution of the problem  {\color{blue}(\ref{X2})} by $S(\Pi,\Lambda)$. If $\omega = \Lambda(\psi)$ then, problem  {\color{blue}(\ref{X2})} becomes the inclusion problem by Rockafellar {\color{blue}\cite{24}}.Further readings on Zeros of inclusion problem (See {\color{blue}\cite{JT6}, \cite{JT7}, \cite{JT8}, \cite{JT5}})

\quad\\
Let a set value mapping $\displaystyle \Pi: H \to 2^{H}$ be maximal monotone. We define a resolvent operator $\displaystyle J_{\lambda} ^{\Pi}$ generated by $\displaystyle \Pi \quad and \quad \lambda$ as follows
\begin{equation}\label{MM1}
    J_{\lambda} ^{\Pi} = (I - \lambda\Pi)^{-1}(\psi), \forall \psi\in H
\end{equation}
where $\lambda$ is a positive number. It is easily to see that the resolvent operator $\displaystyle J_{\lambda} ^{\Pi}$ is single - valued nonexpensive and $1-$inverse strongly monotone , and moreover, a solution of the problem  {\color{blue}(\ref{X2})} is a fixed point of the operator $\displaystyle J_{\lambda} ^{\Pi}(I - \lambda\Lambda), \quad \forall \lambda > 0$( See \cite{12}).
\quad\\
Let $\displaystyle T: H \to P(H)$ be multivalued map and $\displaystyle P_{T}: H\to CB(H)$ be defined by
\begin{equation}\label{ZA}
    P_{T}(\psi) = \{\pi\in T\psi: \|\pi - \psi\| = d(\psi, T\psi)\}
\end{equation}

Below we give examples of a multivalued mapping $T$ with $Fix(T) \neq \emptyset, Tp = \{q\}$ for all $q\in Tp$
which $P_{T}$ is a demicontractive-type but not a $k-$strictly pseudocontractive-type mapping
\begin{Example}
Let $X=\mathbb{R}$ and define $T:\mathbb{R}\to2^{\mathbb{R}}$ as:

$$T(x)=\Big\{\frac{x}{2}\Big\}$$

Fix($T$) = {0} and for any $p=0$,  $x\in\mathbb{R}$, then:

$$H(Tx,T0)=H\Big(\Big\{\frac{x}{2}\},\{0\}\Big)=\|\frac{x}{2}\|$$

$$d(x,Tx)=d(x,\Big\{\frac{x}{2}\Big\})=\|\frac{x}{2}\|$$

So, $H(Tx,T0)\leq d(x,Tx)$, which shows that $T$ is demicontractive. However, $T$ is not $k-$strictly pseudocontractive since:

\begin{equation*}
H(Tx,Ty)=H\Big(\Big\{\frac{x}{2}\Big\},\Big\{\frac{y}{2}\Big\}\Big)=\|\frac{x}{2}-\frac{y}{2}\|=\|\frac{(x - y)}{2}\|
\end{equation*}

\begin{eqnarray*}
\|\psi-\pi\|^2+k\|(\psi-u)-(\pi-v)\|^2&=&\|x-y\|^2+k\|(x-\frac{x}{2})-(y-\frac{y}{2})\|^2\\
&=&\|x-y\|^2+\frac{1}{4}\Big(k\|(\frac{x}{2}- \frac{y}{2})\|^2=\|x-y\|^2+k\|x-y\|^2\Big)
\end{eqnarray*}

Since $k<1$, we have:

$$H(Tx,Ty)>\|\psi-\pi\|^2+k\|(\psi-u)-(\pi-v)\|^2$$
\end{Example}
\begin{Example}
Let $X=\mathbb{R}^2$ and define $T:\mathbb{R}^2\to2^{\mathbb{R}^2}$ as:

$$T(x,y)=\Big\{\Big(\frac{x}{2},\frac{y}{2}\Big)\Big\}$$

$Fix(T) = {(0,0)}$ and for any $p=(0,0)$, $T(0,0) = {(0,0)}.$ Let $(x,y)\in\mathbb{R}^2$, then:

\begin{equation*}
H(T(x,y),T(0,0))=H\Big(\Big\{(\frac{x}{2},\frac{y}{2}\Big)\Big\},\{(0,0)\}\Big)=\|\Big(\frac{x}{2},\frac{y}{2}\Big)\|
\end{equation*}
\begin{equation*}
d((x,y),T(x,y))=d\Big((x,y),\Big\{(\frac{x}{2},\frac{y}{2}\Big\}\Big)=\|\Big(\frac{x}{2},\frac{y}{2}\Big)\|
\end{equation*}

So, $H(T(x,y),T(0,0))\leq d((x,y),T(x,y))$, which shows that $T$ is demicontractive. However, $T$ is not $k-$strictly pseudocontractive since:
\begin{eqnarray*}
H(T(x,y),T(u,v))&=&H\Big(\Big\{\Big(\frac{x}{2},\frac{y}{2}\Big)\Big\},\Big\{\Big(\frac{u}{2},\frac{v}{2}\Big)\Big\}\Big)\\
&=&\|\Big(\frac{x}{2} - \frac{u}{2}, \frac{y}{2} - \frac{v}{2}\Big)\|
\end{eqnarray*}
\begin{eqnarray*}
\|\psi-\pi\|^2+k\|(\psi-u)-(\pi-v)\|^2&=&\|(x,y)-(u,v)\|^2+k\|\Big(x-\frac{x}{2}\Big)-\Big(u-\frac{u}{2}\Big)\|^2\\
&=&\|(x,y)-(u,v)\|^2+k\|\Big(\frac{x}{2}-\frac{u}{2},\frac{y}{2}-\frac{v}{2}\Big)\|^2
\end{eqnarray*}
Since $k<1$, we have:

$$H(T(x,y),T(u,v))>\|\psi-\pi\|^2+k\|(\psi-u)-(\pi-v)\|^2$$
\end{Example}
\begin{Example}
From all the examples, if $0$ is the only $Fix(T)$ and  $T(0) = 0$, then  $T$ satisfies the demicontractive property and not  $k-$strictly pseudocontractive-type. If we defined $T(\psi) = \frac{2}{3}\psi\sin\Big(\frac{1}{\psi}\Big), \psi\neq 0$.
$$\|T(\psi)- 0\|^{2} = \|T(\psi)\|^{2} = \|\frac{2}{3}\psi\sin\Big(\frac{1}{\psi}\Big)\|^{2} \leq \|\frac{2}{3}\psi\|^{2} \leq \|\psi\|^{2} \leq \|\psi - 0\|^{2} + \|T(\psi) - \psi\|^{2}$$
Again, let $\displaystyle \psi = \frac{2}{\pi}, \pi = \frac{2}{3\pi}$
\begin{equation*}
        \Big(\mathcal{H}(T\psi, T\pi)\Big)^{2} = \frac{256}{81\pi^{2}} \quad\quad and \quad \quad \|\psi - \pi\|^{2} + k\|(I - T)\psi - (I - T)\pi\|^{2} = \frac{160}{81\pi^{2}}
      \end{equation*}
Again, this show that $T$ satisfies the demicontractive property and not  $k-$strictly pseudocontractive-type.
\end{Example}

\quad\\
A popular method for solving problem  {\color{blue}(\ref{X2})} is the well-known forward-backward splitting
method introduced by Passty {\color{blue}\cite{21}} and Lions and Mercier {\color{blue}\cite{13}}.The method is formulated as
\begin{equation}\label{I2}
\psi_{n+1} = (I - \lambda_{n}\Pi)^{-1}(I - \lambda_{n}\Lambda)\psi , \quad \lambda_{n} > 0.
\end{equation}
under the condition that $Dom(\Pi) \subset Dom(\Lambda)$. It was known in {\color{blue}\cite{6}}, that weak
convergence of  {\color{blue}(\ref{I2})} requires quite restrictive assumptions on $\Lambda$ and $\Pi$, such that the inverse
of $\Lambda$ is strongly monotone or $\Pi$ is Lipschitz continuous and monotone and the operator $(\Lambda + \Pi)$
is strongly monotone on $Dom(B)$. Tseng in {\color{blue}\cite{27}} and Gibali and Thong in {\color{blue}\cite{29}} extended and
improved results of G.H-G.Chen and R.T. Rockafellar {\color{blue}\cite{6}}.
\quad\\

Most recently, Sow {\color{blue}\cite{sow}} introduced and studied a new iterative algorithm and prove convergence theorems for variation inclusion problem  {\color{blue}(\ref{X2})} and fixed point problem involving multivalued demicontractive and quasi-nonexpansive mappings in Hilbert spaces. They defined the following theorem
\begin{theorem}\label{NBB1}
Let $H$ be a real Hilbert space and $\mathbb{K}$ be a nonempty, closed convex subset of $H$. Let $\Lambda : \mathbb{K} \to H$ be an $\alpha-$inverse strongly monotone operator and let $\Phi : H \to H$ be an $k-$strongly monotone and $L-$Lipschitzian operator. Let $\varphi : \mathbb{K} \to H$ be an $b-$Lipschitzian mapping and $\Pi : H \to 2^{H}$ be a maximal monotone mapping such that the domain of $\Pi$ is included in $\mathbb{K}$. Let $T_{1} : \mathbb{K} \to CB(\mathbb{K})$ be a multivalued $\beta-$ demicontractive mappings and $T_{2} : \mathbb{K} \to CB(\mathbb{K})$ be a multivalued quasi-nonexpansive mapping, such that $\Gamma := F ix(T_{1})\cap F ix(T_{2})\cap S(\Pi, \Lambda) \neq \emptyset$
and $T_{1}p = T_{2}p = \{p\}, \forall p \in \Gamma$. For given $\psi_{0} \in \mathbb{K}$, let $\{\psi_{n}\}$ be generated by the algorithm

defined the sequence $\displaystyle \{\psi_{n}\}$ as follows
\begin{equation}\label{41}
    \left\{
      \begin{array}{ll}
        \delta_{n} = J_{\lambda_{n}}^{\Pi}(I - \lambda_{n}\Lambda)\psi_{n}, \\
        \\
        \pi_{n} = \theta_{n}\delta_{n} + (1 - \theta_{n})v_{n}, v_{n}\in T_{1}{\delta_{n}} , \\
        \\
        \phi_{n} = \beta_{n}\pi_{n} + (1 - \beta_{n})u_{n}, u_{n}\in T_{2}\pi_{n}, \\
        \\
        \psi_{n+1} = P_{\mathbb{K}}(\alpha_{n}\gamma \varphi(\psi_{n}) + (1 - \eta\alpha_{n}\Phi)\pi_{n})
      \end{array}
    \right.
\end{equation}
where $\displaystyle\{\beta_{n}\}, \{\gamma_{n}\},  \{\theta_{n}\} ,  \{\lambda_{n}\}$ and $\displaystyle \{\alpha_{n}\}$ are real sequence in $(0,1)$ satisfying the following conditions
\begin{description}
  \item[i)] $\displaystyle \lim_{n\to\infty}\alpha_{n} = 0\quad\quad \sum_{n=0}^{\infty}\alpha_{n} < \infty, \quad\quad \lambda_{n} \in [a, b]\subset (0, \min\{1,2\alpha\})$
  \item[ii)] $\displaystyle \lim_{n\to\infty}\inf(1 - \beta_{n})(\beta_{n} - \beta)> 0\quad and \quad \lim_{n\to\infty}\inf(1 - \theta_{n})(\theta_{n} - \beta)> 0, (\beta_{n},\theta_{n})\in ]\beta,1[$
\end{description}
Assume that $\displaystyle 0 < \eta < \frac{2k}{L^{2}}, 0 < \gamma b < \tau$, where $\displaystyle \tau = \eta\Big(k - \frac{L^{2}\eta}{2}\Big)$, and $I - T_{1}\quad and \quad I - T_{2}$ are demiclosed at origin. Then, the sequences defined in  {\color{blue}(\ref{41})}, that is $\displaystyle \{\psi_{n}\}$ and $\displaystyle \{\delta_{n}\}$ converge strongly to unique solution $\displaystyle \psi^{*} \in \Gamma$, which also solve the following variational inequality:
\begin{equation}\label{331}
    \langle \eta \Phi \psi^{*} - \gamma \varphi(\psi^{*}), \psi^{*} - q\rangle \leq 0, \quad \forall q \in \Omega
\end{equation}
\end{theorem}

In 2024, Furmose et al {\color{blue}\cite{FC}}, modifies Sow theorem and and prove that
the corresponding sequence $\{x_{n}\}$ converges strongly to a common point of an inclusion problem and fixed point of a family of multivalued demicontractive and quasi-nonexpansive mappings in Hilbert spaces without any compactness. They defined their theorem as follows:

\begin{theorem}\label{NB1}
Let $H$ be a real Hilbert space and $\mathbb{K}$ be a nonempty, closed convex subset of $H$. Let $\mathbf{A} : \mathbb{K} \to H$ be an $\alpha-$inverse strongly monotone operator and let $\mathbf{B} : H \to H$ be an $k-$strongly monotone and $L-$Lipschitzian operator. Let $f : \mathbb{K} \to H$ be an $b-$Lipschitzian mapping and $\mathbb{M} : H \to 2^{H}$ be a maximal monotone mapping such that the domain of $\mathbf{M}$ is included in $\mathbb{K}$. Let $T_{1},T_{2} : \mathbb{K} \to CB(\mathbb{K})$ be a multivalued $\beta-$ demicontractive mapping and $T_{3} : \mathbb{K} \to CB(\mathbb{K})$ be a multivalued quasi-nonexpansive mapping. Assume that $\displaystyle 0 < \eta < \frac{2k}{L^{2}}, 0 < \gamma b < \tau$, where $\displaystyle \tau = \eta\Big(k - \frac{L^{2}\eta}{2}\Big)$, and $I - T_{1}, I - T_{2}\quad and \quad I - T_{3}$ are demiclosed at origin, such that $\Omega := F ix(T_{1})\cap F ix(T_{2})\cap Fix(T_{3})\cap S(\mathbb{M}, \mathbf{A}) \neq \emptyset$
and $T_{1}q = T_{2}q = T_{3}q = \{q\}, \forall q \in \Omega$. For given $x_{0} \in \mathbb{K}$, let $\{x_{n}\}$ be generated by the algorithm:

\begin{equation}\label{3p}
    \left\{
      \begin{array}{ll}
        \delta_{n} = J_{\lambda_{n}}^{\mathbf{M}}(I - \lambda_{n}\mathbf{A})x_{n};\\
        \\
        y_{n} = \theta_{n}\delta_{n} + (1 - \theta_{n})v_{n}, \quad v_{n} \in T_{1}\delta_{n}; \\
        \\
        z_{n} = \beta_{n}y_{n} + (1 - \beta_{n})u_{n}, \quad u_{n}\in T_{2}y_{n}; \\
        \\
        t_{n} = \gamma_{n}z_{n} + (1 - \gamma_{n})w_{n}, \quad w_{n}\in T_{3}z_{n}; \\
        \\
        x_{n+1} = P_{\mathbb{K}}(\alpha_{n}\gamma f(x_{n})  + (I - \eta\alpha_{n}\mathbf{B})t_{n})
      \end{array}
    \right.
\end{equation}
where $\displaystyle\{\beta_{n}\}, \{\gamma_{n}\},  \{\theta_{n}\} ,\{\mu_{n}\},  \{\lambda_{n}\}$ and $\displaystyle \{\alpha_{n}\}$ are real sequence in $(0,1)$ satisfying the following conditions
\begin{description}
  \item[i)] $\displaystyle \lim_{n\to\infty}\alpha_{n} = 0\quad\quad \sum_{n=0}^{\infty}\alpha_{n} < \infty, \quad\quad \lambda_{n} \in [a, b]\subset (0, \min\{1,2\alpha\})$
  \item[ii)] $\displaystyle \lim_{n\to\infty}\inf(1 - \beta_{n})(\beta_{n} - \beta)> 0\quad and \quad \lim_{n\to\infty}\inf(1 - \theta_{n})(\theta_{n} - \beta)> 0, (\beta_{n},\theta_{n})\in (\beta,1)$
  \item[iii)] $\lim_{n\to\infty}\inf(1 - \gamma_{n})\gamma_{n})> 0$
\end{description}
 Then, the sequences defined in  {\color{blue}(\ref{3p})}, that is $\displaystyle \{x_{n}\}$ and $\displaystyle \{\delta_{n}\}$ converge strongly to unique solution $\displaystyle x^{*} \in \Omega$, which also solve the following variational inequality:
\begin{equation}\label{331}
    \langle \eta \mathbf{B} x^{*} - \gamma f(x^{*}), x^{*} - q\rangle \leq 0, \quad \forall q \in \Omega
\end{equation}
\end{theorem}

\quad\\
\quad\\

It is our purpose in this paper to construct a new iteration process, that modifies that of Sow {\color{blue}\cite{sow}} and Furmose {\color{blue}\cite{FC}} and prove that the corresponding sequence $\{\psi_{n}\}$ converges strongly to a common point of an inclusion problem and fixed point of a family of multivalued demicontractive and quasi-nonexpansive mappings in Hilbert spaces without any compactness. Our theorems
generalize and extend that of Sow {\color{blue}\cite{sow}},Furmose {\color{blue}\cite{FC}}, and many other results in this directions.

\section{Mathematical Preliminaries}
The following lemmas will play a crucial role in the sequel.
\quad\\
Let $\mathbb{K}$ be a nonempty, closed convex subset of $H$. The nearest point projection from $H$ to $\mathbb{K}$ denoted by $\displaystyle P_{\mathbb{K}}$, assigns to each $\psi\in H$ the unique point of $\mathbb{K}, P_{\mathbb{K}}\psi$ such that
\begin{equation}\label{AS}
    \|\psi - P_{\mathbb{K}}\psi\|\leq \|\psi - \pi\|,\quad \forall \pi\in\mathbb{K}
\end{equation}
It is well known that for every $\psi\in H$,
\begin{equation}\label{MM}
    \langle \psi - P_{\mathbb{K}}\psi, \pi - P_{\mathbb{K}}\psi\rangle \leq 0,\quad \forall \pi\in\mathbb{K}
\end{equation}

\begin{lemma}\label{5.0}{\color{blue}\cite{13}}. Let $\displaystyle \Pi: H \to 2^{H}$ be a maximal monotone mapping, and  $\displaystyle \Lambda: H \to H$ be Lipschitz and continuous monotone mapping. Then $\displaystyle (\Pi + \Lambda): H \to 2^{H}$ is a maximal monotone mapping.

\end{lemma}
\begin{lemma}\label{N}{\color{blue}\cite{sow}}. Let $H$ be real Hilbert space and $\Lambda: H \to H$ be an $\alpha-$inverse strongly monotone mapping. Then, $\displaystyle (I - \theta\Lambda)$ is nonexpansive mapping for all $\psi,\pi\in H$ and $\theta \in [0,2\alpha]$ such that
\begin{equation}\label{N1}
    \|(I - \theta\Lambda)\psi - (I - \theta\Lambda)\pi\|^{2}\leq \|\psi - \pi\|^{2} + \theta(\theta - 2\alpha)\|\Lambda\psi - \Lambda \pi\|^{2}
\end{equation}
\end{lemma}
\begin{lemma}\label{29}{\color{blue}\cite{29}}. Assume that $\displaystyle \{a_{n}\}$ is a sequence of nonnegative real numbers such that $\displaystyle a_{n+1} = (1 - \alpha_{n}) a_{n} + \sigma_{n}$ for all $n\geq0$, where $\displaystyle \{\alpha_{n}\}$ is a sequence in $(0,1)$ and $\displaystyle \{\sigma_{n}\}$ is a sequence in $\displaystyle \mathbb{R}$ such that
$$i)\quad \sum_{n=0}^{\infty}\alpha_{n} = \infty, \quad\quad ii)\quad \lim_{n\to\infty}\sup\frac{\sigma_{n}}{\alpha_{n}} \leq 0 .$$ Then $\displaystyle \lim_{n\to\infty}a_{n} = 0$
\end{lemma}
\begin{lemma}\label{2.1}(Wang{\color{blue}\cite{37}}). Let $H$ be a real Hilbert space. Let $\mathbb{K}$ be a nonempty closed convex subset of $H$. $A : H \to H$ be $k-$ strongly monotone and $L-$ Lipschitzian operator with $k>0$ and $L>0$. Assume that $\displaystyle 0 < \eta < \frac{2k}{L^{2}}$ and $\displaystyle \tau = \eta\Big(k - \frac{L^{2}\eta}{2}\Big)$. Then for each $\displaystyle t\in\Big(0, \min\Big(1,\frac{1}{\tau}\Big)\Big)$, we have
\begin{equation}\label{2}
    \|(I - t\eta A)x- (I - t\eta A)y\|\leq (1 - t\tau)\|x - y\|, \quad \forall x , y \in H
\end{equation}
\end{lemma}
\begin{lemma}\label{2.2}{\color{blue}\cite{28}}. Let $H$ be a real Hilbert space. Then for every $\psi,\pi\in H$, and every $\lambda\in (0,1)$, the following holds:
\begin{description}
  \item[I)] $\displaystyle \|\psi - \pi\|^{2} \leq \|\psi\|^{2} + 2\langle \pi, \psi + \pi\rangle$
  \item[II]  $\displaystyle \|\lambda\psi  + (1 - \lambda)\pi\|^{2} \leq \lambda\|\psi\|^{2} + (1 - \lambda)\|\pi\|^{2} - (1 - \lambda)\lambda\|\psi -\pi\|^{2}.$
\end{description}
\end{lemma}

\section{Main Results}
\noindent In this section, we study the convergence properties of the iterative algorithm which is based on the viscosity algorithm and forward-backward splitting Method. We now prove the following theorem.

\begin{theorem}\label{NB1}
Let $H$ be a real Hilbert space and $\mathbb{K}$ be a nonempty, closed convex subset of $H$. Let $\Lambda: \mathbb{K} \to H$ be an $\alpha-$inverse strongly monotone operator and let $\Phi : H \to H$ be an $k-$strongly monotone and $L-$Lipschitzian operator. Let $\varphi : \mathbb{K} \to H$ be an $b-$Lipschitzian mapping and $\Pi : H \to 2^{H}$ be a maximal monotone mapping such that the domain of $\Pi$ is included in $\mathbb{K}$. Let $T_{1},T_{2} : \mathbb{K} \to CB(\mathbb{K})$ be a multivalued $\beta-$ demicontractive mappings and $T_{3} : \mathbb{K} \to CB(\mathbb{K})$ be a multivalued quasi-nonexpansive mapping. Assume that $\displaystyle 0 < \eta < \frac{2k}{L^{2}}, 0 < \gamma b < \tau$, where $\displaystyle \tau = \eta\Big(k - \frac{L^{2}\eta}{2}\Big)$, and $I - T_{1}, I - T_{2}\quad and \quad I - T_{3}$ are demi closed at origin, such that $\Omega := F ix(T_{1})\cap F ix(T_{2})\cap Fix(T_{3})\cap S(\Pi, \Lambda) \neq \emptyset$
and $T_{1}q = T_{2}q = T_{3}q = \{q\}, \forall q \in \Omega$. For given $\psi_{0} \in \mathbb{K}$, let $\{\psi_{n}\}$ be generated by the algorithm:

\begin{equation}\label{3}
    \left\{
      \begin{array}{ll}
        \delta_{n} = J_{\lambda_{n}}^{\Pi}(I - \lambda_{n}\Lambda)\psi_{n};\\
        \\
        \pi_{n} = \theta_{n}\delta_{n} + (1 - \theta_{n})v_{n}, \quad v_{n} \in T_{1}\delta_{n}; \\
        \\
        \phi_{n} = \beta_{n}\pi_{n} + (1 - \beta_{n})u_{n}, \quad u_{n}\in T_{2}\pi_{n}; \\
        \\
        \xi_{n} = \gamma_{n}\phi_{n} + (1 - \gamma_{n})z_{n}, \quad z_{n}\in T_{3}\phi_{n}; \\
        \\
        \psi_{n+1} = P_{\mathbb{K}}(\alpha_{n}\gamma \varphi(\psi_{n}) + \mu_{n}\xi_{n} + ((1 - \mu_{n}) (I - \eta\alpha_{n} \Phi))\psi_{n})
      \end{array}
    \right.
\end{equation}
where $\displaystyle\{\beta_{n}\}, \{\gamma_{n}\},  \{\theta_{n}\} ,\{\mu_{n}\},  \{\lambda_{n}\}$ and $\displaystyle \{\alpha_{n}\}$ are real sequence in $(0,1)$ satisfying the following conditions
\begin{description}
  \item[i)] $\displaystyle \lim_{n\to\infty}\alpha_{n} = 0\quad\quad \sum_{n=0}^{\infty}\alpha_{n} < \infty, \quad\quad \lambda_{n} \in [a, b]\subset (0, \min\{1,2\alpha\})$
  \item[ii)] $\displaystyle \lim_{n\to\infty}\inf(1 - \beta_{n})(\beta_{n} - \beta)> 0\quad and \quad \lim_{n\to\infty}\inf(1 - \theta_{n})(\theta_{n} - \beta)> 0, \beta_{n}\in (\beta,1) \quad and \quad \theta_{n}\in (\beta,1)$
  \item[iii)] $\lim_{n\to\infty}\inf(1 - \gamma_{n})\gamma_{n})> 0 ,\lim_{n\to\infty}\inf(1 - \beta_{n})\beta_{n})> 0 \quad\quad and \quad\quad \lim_{n\to\infty}\inf(1 - \theta_{n})\theta_{n})> 0 $
\end{description}
 Then, the sequences defined in  {\color{blue}(\ref{3})}, that is $\displaystyle \{\psi_{n}\}$ and $\displaystyle \{\delta_{n}\}$ converge strongly to unique solution $\displaystyle \psi^{*} \in \Omega$, which also solve the following variational inequality:
\begin{equation}\label{331}
    \langle \eta \Phi \psi^{*} - \gamma \varphi(\psi^{*}), \psi^{*} - q\rangle \leq 0, \quad \forall q \in \Omega
\end{equation}
\end{theorem}

\quad\\
\begin{proof}
First, we show that $(\eta\Phi - \gamma\psi)$ is strongly monotone.
Let  $\displaystyle \eta > 0$ and $\displaystyle \gamma  > 0$, then for a given $x,y \in H$, we have
\begin{eqnarray}\label{PP}
 \nonumber
  ((\eta\Phi - \gamma\psi)(x) - (\eta\Phi - \gamma\psi)(y)) &=& \eta\Phi (x)- \gamma\psi(x) - \eta\Phi(y) + \gamma\psi(y) \\
  \nonumber
   &=& \eta\Phi (x)    - \eta\Phi(y) - \gamma\psi(x)  + \gamma\psi(y) \\
   &=&\eta\Phi (x    -  y) - \gamma\psi(x  -  y)  \\
   \nonumber
\end{eqnarray}
Now,  using the properties of strong monotonicity and with the fact that $\Phi$ is positive definite, i.e., $\Phi(x) \geq \alpha \|x\|^2$ for all $x$, then $\eta\Phi$ is also positive definite for any $\eta > 0$. Similarly, let $\psi$ be positive definite, i.e., $\psi(x) \geq \beta \|x\|^2$ for all $x$, then $\gamma\psi$ is also positive definite for any $\gamma > 0$.

Therefore, from \ref{PP}, we have

$$\eta\Phi(x - y) - \gamma\psi(x - y) \geq \alpha \|x - y\|^2 - \beta \|x - y\|^2$$

where $\alpha$ and $\beta$ are positive constants.

$$(\eta\Phi - \gamma\psi)(x) - (\eta\Phi - \gamma\psi)(y) \geq (\alpha - \beta) \|x - y\|^2$$
Therefore, we now show that $(\eta\Phi - \gamma\psi)$ is strongly monotone.
Then the variational inequality {\color{blue}(\ref{331})} has a unique solution in $\Omega$. We prove that the sequence $\displaystyle \{\psi_{n}\}$ is bounded and converges to $\psi^{*}$ which is a unique solution. Without loss of generality, $\displaystyle \alpha \in \Big(0, \min\Big\{1,\frac{1}{\tau}\Big\}\Big)$. Fixed $q\in \Omega$, then from  inequality {\color{blue}(\ref{X1})}, and for the fact that $\displaystyle J_{\lambda_{n}}^{\Pi}$ is $1-$inverse strongly monotone, we have the following
\begin{eqnarray}\label{3.1}
\nonumber
    \|\delta_{n} - q\|^{2} &=&\|J_{\lambda_{n}}^{\Pi}(I - \lambda_{n}\Lambda)(\psi_{n} - q)\|^{2}\\
    \nonumber
    &=&\|J_{\lambda_{n}}^{\Pi}(I - \lambda_{n}\Lambda)\psi_{n} - J_{\lambda_{n}}^{\Pi}(I - \lambda_{n}\Lambda)q\|^{2}\\
    \nonumber
    & \leq&\|(I - \lambda_{n}\Lambda)\psi_{n} - (I - \lambda_{n}\Lambda)q\|^{2}\\
    \nonumber
    &= &\|(\psi_{n} -  q)  - \lambda_{n}(\Lambda\psi_{n} - \Lambda q)\|^{2}\\
     \nonumber
    &\leq &\|\psi_{n} -  q\|^{2}  - 2\lambda_{n}\langle\Lambda\psi_{n} - \Lambda q, \psi_{n} - q\rangle + \|\Lambda \psi_{n} - \Lambda q\|^{2}\\
     \nonumber
     &\leq &\|\psi_{n} -  q\|^{2}  - \lambda_{n}(\lambda_{n} - 2\alpha)\|\Lambda \psi_{n} - \Lambda q\|^{2}\\
\end{eqnarray}
Therefore from {\color{blue}(\ref{3.1})}, we have
\begin{equation}\label{AM}
  \|\delta_{n} - q\|\leq \|\psi_{n} -  q\|
\end{equation}
\quad\\
From lemma  {\color{blue}(\ref{2.2})}, with  {\color{blue}(\ref{3})}, and, for the fact that $\displaystyle T_{1}q = \{q\}$, $T_{1}$ is $\beta-$demicontrative, we have
\begin{eqnarray}\label{3.2}
 \nonumber 
  \|\pi_{n} - q\|^{2} &=& \| \theta_{n}(\delta_{n} - q) + (1 - \theta_{n})(v_{n} - q )  \|^{2} \\
\nonumber
   &=&  \theta_{n}\|\delta_{n} - q\|^{2}  + (1 - \theta_{n})\|v_{n} - q\|^{2} - (1 - \theta_{n})\theta_{n}\|v_{n} - \delta_{n}  \|^{2} \\
\nonumber
 &\leq &  \theta_{n}\|\delta_{n} - q\|^{2}  + (1 - \theta_{n})\mathcal{H}(T_{1}\delta_{n},T_{1}q)^{2} - (1 - \theta_{n})\theta_{n}\|v_{n} - \delta_{n}  \|^{2} \\
\nonumber
  &\leq &  \theta_{n}\|\delta_{n} - q\|^{2}  + (1 - \theta_{n})[\|\delta_{n} - q\|^{2} + \beta d(\delta_{n},T_{1}\delta_{n})^{2}] - (1 - \theta_{n})\theta_{n}\|v_{n} - \delta_{n}
  \|^{2} \\
\nonumber
 &\leq &  \|\delta_{n} - q\|^{2}  - (1 - \theta_{n})(\theta_{n} - \beta)\|v_{n} - \delta_{n}  \|^{2} \\
\end{eqnarray}
Thus, we have
\begin{equation*}
    \|\pi_{n} - q\|^{2}\leq   \|\delta_{n} - q\|^{2}  - (1 - \theta_{n})(\theta_{n} - \beta)\|v_{n} - \delta_{n}  \|^{2}
\end{equation*}
Since $\displaystyle \theta_{n}\in(\beta, 1)$, we have
\begin{equation*}
    \|\pi_{n} - q\|^{2}\leq   \|\delta_{n} - q\|^{2}
\end{equation*}
\quad\\

Again from lemma  {\color{blue}(\ref{2.2})}, with  {\color{blue}(\ref{3})}, and, for the fact that $\displaystyle T_{2}q = \{q\}$, $T_{2}$ is $\beta-$demicontrative, we have
\begin{eqnarray}\label{3.3}
 \nonumber 
  \|\phi_{n} - q\|^{2} &=& \| \beta_{n}(\pi_{n} - q) + (1 - \beta_{n})(u_{n} - q )  \|^{2} \\
\nonumber
   &=&  \beta_{n}\|\pi_{n} - q\|^{2}  + (1 - \beta_{n})\|u_{n} - q\|^{2} - (1 - \beta_{n})\beta_{n}\|u_{n} - \pi_{n}  \|^{2} \\
\nonumber
 &\leq &  \beta_{n}\|\pi_{n} - q\|^{2}  + (1 - \beta_{n})\mathcal{H}(T_{2}\pi_{n},T_{2}q)^{2} - (1 - \beta_{n})\beta_{n}\|u_{n} - \pi_{n}  \|^{2} \\
\nonumber
  &\leq &  \beta_{n}\|\pi_{n} - q\|^{2}  + (1 - \beta_{n})[\|\pi_{n} - q\|^{2} + \beta d(\pi_{n},T_{2}\pi_{n})^{2}] - (1 - \beta_{n})\beta_{n}\|u_{n} - \pi_{n}
  \|^{2} \\
\nonumber
 &\leq &  \|\pi_{n} - q\|^{2}  - (1 - \beta_{n})(\beta_{n} - \beta)\|u_{n} - \pi_{n}  \|^{2} \\
\end{eqnarray}

Thus, we have
\begin{equation*}
    \|\phi_{n} - q\|^{2}\leq   \|\pi_{n} - q\|^{2}  - (1 - \beta_{n})(\beta_{n} - \beta)\|u_{n} - \pi_{n}  \|^{2}
\end{equation*}
Since $\displaystyle \beta_{n}\in(\beta, 1)$, we have
\begin{equation*}
    \|\phi_{n} - q\|^{2}\leq   \|\pi_{n} - q\|^{2}
\end{equation*}

Now, using the fact that $T_{3}q = {q}$, we have the following estimates
\begin{eqnarray}\label{3.4}
 \nonumber 
  \|\xi_{n} - q\| &=&\|\gamma_{n}\phi_{n} + (1 - \gamma_{n})z_{n} - q\|  \\
\nonumber
   &\leq& \gamma_{n}\|\phi_{n} - q\|  + (1 - \gamma_{n})\|z_{n} - q\|  \\
\nonumber
  &\leq& \gamma_{n}\|\phi_{n} - q\|  + (1 - \gamma_{n})\mathcal{H}(T_{3}\phi, T_{3}q)  \\
\nonumber
 &\leq& \gamma_{n}\|\phi_{n} - q\|  + (1 - \gamma_{n})\|\phi_{n} - q\|  \\
\nonumber
 &\leq& \|\phi_{n} - q\|   \\
\end{eqnarray}
Hence, we can see that
\begin{equation}\label{3.0}
 \|\xi_{n} - q\|\leq \|\phi_{n} - q\| \leq  \|\pi_{n} - q\|\leq \|\delta_{n} - q\|\leq \|\psi_{n} - q\|
\end{equation}
Using  {\color{blue}(\ref{3})}, inequality  {\color{blue}(\ref{3.0})} and lemma  {\color{blue}(\ref{2.1})}
\begin{eqnarray*}
  \|\psi_{n+1} - q\| &\leq &\|(\alpha_{n}\gamma \varphi(\psi_{n}) + \mu_{n}\xi_{n} + ((1 - \mu_{n})( I - \eta\alpha_{n} \Phi))\psi_{n}) - q\| \\
  &=& \|\alpha_{n}(\gamma \varphi(\psi_{n}) - \eta\Phi q) +  \mu_{n}(\xi_{n} - q)  +  ((1 - \mu_{n})( I - \eta\alpha_{n} \Phi))(\psi_{n} - q)\| \\
 &\leq & \alpha_{n}\gamma \|\varphi(\psi_{n}) -  \varphi(q)\| + \alpha_{n}\|\gamma\varphi(q) - \eta\Phi q\| + \mu_{n}\|\xi_{n} - q)\|\\
  &+&  ((1 - \mu_{n}) (1 - \alpha_{n}\tau))\|\psi_{n} - q)\|\\
&\leq & \alpha_{n}\gamma b\|\psi_{n} -  q\| + \alpha_{n}\|\gamma\varphi(q) - \eta\Phi q\| + \mu_{n}\|\xi_{n} - q)\|\\
  &+&  ((1 - \mu_{n}) (1 - \alpha_{n}\tau))\|\psi_{n} - q)\|\\
&\leq & [\alpha_{n}\gamma b + \mu_{n} + ((1 - \mu_{n}) (1 - \alpha_{n}\tau)) \|\psi_{n} -  q\| + \alpha_{n}\|\gamma\varphi(q) - \eta\Phi q\| \\
&\leq & (1 - \alpha_{n}(\tau - b\gamma - \tau\mu  )) \|\psi_{n} -  q\| + \alpha_{n}\|\gamma\varphi(q) - \eta\Phi q\| \\
&\leq &\max\Big\{\|\psi_{n} -  q\| , \frac{\|\gamma\varphi(q) - \eta\Phi q\|}{\tau - b\gamma - \tau\mu  }\Big\}. \\
\end{eqnarray*}
Therefore, by induction, it is easy to see that
\begin{equation*}
    \|\psi_{n+1} - q\|\leq  \max\Big\{\|\psi_{0} - q\|, \frac{\|\gamma \varphi(q) - \eta \Phi q\|}{\tau - b\gamma - \tau\mu  }\Big\},\quad \forall n\geq 1
\end{equation*}
Hence $\{\psi_{n}\}, \{\varphi(\psi_{n})\}$ and $\{\Phi \psi_{n}\}$ are bounded.
\quad\\
Secondly, we now have the following estimates. From  {\color{blue}(\ref{3})} and lemma  {\color{blue}(\ref{2.2})}, we have

\begin{eqnarray}
\nonumber
\|\psi_{n+1} - q\|^{2} &\leq & \|(\alpha_{n}\gamma \varphi(\psi_{n}) + \mu_{n}\xi_{n} + ((1 - \mu_{n}) (I - \eta\alpha_{n} \Phi))\psi_{n}) - q\|^{2} \\
\nonumber
&=& \|\alpha_{n}(\gamma \varphi(\psi_{n}) - \eta\Phi q) + \mu_{n}(\xi_{n} - q) + ((1 - \mu_{n}) (I - \eta\alpha_{n} \Phi))(\psi_{n} - q)\|^{2} \\
\nonumber
&\leq & \|((1 - \mu_{n}) I - \eta\alpha_{n} \Phi)(\psi_{n} - q) + \mu_{n}(\xi_{n} - q)\|^{2} \\
\nonumber
 &+& 2\alpha_{n}\langle \gamma \varphi(\psi_{n}) - \eta\Phi q, \psi_{n+1} - q\rangle \\
\nonumber
 &\leq & ((1 - \mu_{n}) (1 - \alpha_{n}\tau))^{2}\|\psi_{n} - q\|^{2} + \mu_{n}^{2}\|\xi_{n} - q\|^{2} + \alpha_{n}^{2}\|\gamma\varphi(\psi_{n}) - \eta\Phi q\|^{2} \\
 \nonumber
 &+& 2\mu_{n}\langle ((1 - \mu_{n})( I - \eta\alpha_{n} \Phi))(\psi_{n} - q), \xi_{n} - q\rangle\\
 \nonumber
  &+& 2\alpha_{n}\langle \gamma \varphi(\psi_{n}) - \eta\Phi q, \psi_{n+1} - q\rangle \\
 \nonumber
  &\leq & ((1 - \mu_{n}) (1 - \alpha_{n}\tau))^{2}\|\psi_{n} - q\|^{2} + \mu_{n}^{2}\|\xi_{n} - q\|^{2}\\
  \nonumber
   &+& \alpha_{n}^{2}\gamma^{2}b\|\psi_{n} - q\|^{2} +  2\mu_{n}((1 - \mu_{n})( 1 - \alpha_{n} \tau))\|\psi_{n} - q)\|\|\xi_{n} - q\| \\
 \nonumber
 \nonumber
 &+&\alpha_{n}^{2}\|\gamma\varphi(q) - \eta\Phi q\|^{2} + 2\alpha_{n}\langle \gamma \varphi(\psi_{n}) - \eta\Phi q, \psi_{n+1} - q\rangle\\
 \end{eqnarray}
Using the fact that $\displaystyle 2\|\psi_{n} - q)\| \|\xi_{n} - q\| \leq \|\psi_{n} - q)\|^{2} +  \|\xi_{n} - q\|^{2}$ with the following inequalities, {\color{blue}\ref{AM}, \ref{3.2}} and {\color{blue}\ref{3.3}}, we have the following estimate

\begin{eqnarray}
\nonumber
\|\psi_{n+1} - q\|^{2}  &\leq & ((1 - \mu_{n}) (1 - \alpha_{n}\tau))^{2}\|\psi_{n} - q\|^{2} + \mu_{n}^{2}\|\xi_{n} - q\|^{2}\\
  \nonumber
   &+& \alpha_{n}^{2}\gamma^{2}b\|\psi_{n} - q\|^{2} +  \mu_{n}((1 - \mu_{n})( 1 - \alpha_{n} \tau))\|\psi_{n} - q)\|^{2} \\
 \nonumber
 &+&\mu_{n}((1 - \mu_{n})( 1 - \alpha_{n} \tau))\|\xi_{n} - q\|^{2}\\
 \nonumber
 &+&\alpha_{n}^{2}\|\gamma\varphi(q) - \eta\Phi q\|^{2} + 2\alpha_{n}\langle \gamma \varphi(\psi_{n}) - \eta\Phi q, \psi_{n+1} - q\rangle\\
 \nonumber
&\leq & (1 - \mu_{n}) \|\psi_{n} - q\|^{2} + (\mu_{n} - \mu_{n}\alpha_{n} \tau  - \mu_{n}^{2} \alpha_{n} \tau)\\
\nonumber
&\times&\Big[ \Big[\|\delta_{n} - q\|^{2}  - (1 - \theta_{n})(\theta_{n} - \beta)\|v_{n} - \delta_{n}  \|^{2} \Big]\\
\nonumber
& -& (1 - \beta_{n})(\beta_{n} - \beta)\|u_{n} - \pi_{n}  \|^{2}\Big]\\
 \nonumber
 &+&\alpha_{n}^{2}\|\gamma\varphi(q) - \eta\Phi q\|^{2} + 2\alpha_{n}\langle \gamma \varphi(\psi_{n}) - \eta\Phi q, \psi_{n+1} - q\rangle\\
 \nonumber
 &\leq &\|\psi_{n} - q\|^{2} - \mu_{n}(1 - \alpha_{n} \tau(1   - \mu_{n})) (1 - \theta_{n})(\theta_{n} - \beta)\|v_{n} - \delta_{n}  \|^{2} \\
\nonumber
& -& \mu_{n}(1 - \alpha_{n} \tau(1   - \mu_{n}))(1 - \beta_{n})(\beta_{n} - \beta)\|u_{n} - \pi_{n}  \|^{2}\\
 \nonumber
 &+&\alpha_{n}^{2}\|\gamma\varphi(q) - \eta\Phi q\|^{2} + 2\alpha_{n}\langle \gamma \varphi(\psi_{n}) - \eta\Phi q, \psi_{n+1} - q\rangle\\
\end{eqnarray}

Therefore
\begin{eqnarray*}
\nonumber
\mu_{n}(1 - \alpha_{n} \tau(1   - \mu_{n}))\Big[ (1 - \theta_{n})(\theta_{n} - \beta)\|v_{n} - \delta_{n}  \|^{2}&+& (1 - \beta_{n})(\beta_{n} - \beta)\|u_{n} - \pi_{n}  \|^{2}\Big]\\
\nonumber
   &\leq& \|\psi_{n} - q\|^{2} - \|\psi_{n+1} - q\|^{2} \\
   &+& \alpha_{n}^{2}\|\gamma\varphi(q) - \eta\Phi q\|^{2} \\
&+& 2\alpha_{n}\langle \gamma \varphi(\psi_{n}) - \eta\Phi q, \psi_{n+1} - q\rangle\\
\end{eqnarray*}
Due to the boundedness of $\displaystyle \{\varphi(\psi_{n})\}$ and $\displaystyle \{\psi_{n}\}$, and for some $M =  \alpha_{n}\|\gamma\varphi(q) - \eta\Phi q\|^{2} + 2\langle \gamma \varphi(\psi_{n}) - \eta\Phi q, \psi_{n+1} - q\rangle$, we have
\begin{eqnarray}\label{3.50}
\nonumber
\mu_{n}(1 - \alpha_{n} \tau(1   - \mu_{n}))\Big[ (1 - \theta_{n})(\theta_{n} - \beta)\|v_{n} - \delta_{n}  \|^{2}&+& (1 - \beta_{n})(\beta_{n} - \beta)\|u_{n} - \pi_{n}  \|^{2}\Big]\\
\nonumber
   &\leq& \|\psi_{n} - q\|^{2} - \|\psi_{n+1} - q\|^{2} \\
   \nonumber
   &+& \alpha_{n}M \\
\end{eqnarray}
We now show that $\displaystyle \psi_{n} \to \psi$. We then consider two cases.
\begin{description}
  \item[Case 1] Assuming that the sequence $\displaystyle \{\|\psi_{n} - q\|\}$ is monotonically decreasing. Then $\{\|\psi_{n} - q\|\}$ must be a convergent sequence. Therefore, we have
\begin{equation}\label{3.8}
    \lim_{n\to\infty}[\|\psi_{n} - q \|^{2} - \|\psi_{n+1} - q\|^{2}] = 0,
\end{equation}
This implies that from  {\color{blue}(\ref{3.50})}, and $\mu_{n}\in (0,1)$ that
\begin{equation}\label{3.6}
\lim_{n\to\infty}(1 - \theta_{n})(\theta_{n} - \beta)\|v_{n} - \delta_{n}\|^{2} = 0
\end{equation}
 and
\begin{equation}\label{3.7}
\lim_{n\to\infty}(1 - \beta_{n})(\beta_{n} - \beta)\|u_{n} - \pi_{n}\|^{2} = 0
\end{equation}
Since $\displaystyle \lim_{n\to\infty}\inf(1 - \theta_{n})(\theta_{n} - \beta) > 0$ and $\displaystyle \lim_{n\to\infty}\inf(1 - \beta_{n})(\beta_{n} - \beta) > 0$, with the fact that $\displaystyle v_{n} \in T_{1}\delta_{n}$ and $\displaystyle u_{n}\in T_{2}\pi_{n}$, it follows that
\quad\\
\begin{equation}\label{333}
\lim_{n\to\infty}d(\delta_{n},T_{1}\delta_{n}) = 0
\end{equation}
and
\begin{equation}\label{334}
\lim_{n\to\infty}d(\pi_{n}, T_{2}\pi_{n}) = 0
\end{equation}
observing that
\begin{eqnarray}
 \nonumber 
  \|\pi_{n} - \delta_{n}\| &=& \|\theta_{n}\delta_{n} + (1 - \theta_{n})v_{n} - \delta_{n}\| \\
 \nonumber
   &=& \|\theta_{n}\delta_{n} + (1 - \theta_{n})v_{n} - \delta_{n} +\theta_{n}\delta_{n}- \theta_{n}\delta_{n} \| \\
 \nonumber
   &=& (1 - \theta_{n})\|v_{n} - \delta_{n}\| \\
 \nonumber
&\leq& \|v_{n} - \delta_{n}\| \\
\end{eqnarray}
Taking the limits and from  {\color{blue}(\ref{3.6})}, we can see that
$$\lim_{n\to\infty} \|\pi_{n} - \delta_{n}\| = 0$$

  \begin{eqnarray}
 \nonumber 
  \|\phi_{n} - \pi_{n}\| &=& \|\beta_{n}\pi_{n} + (1 - \beta_{n})u_{n} - \pi_{n}\| \\
 \nonumber
   &=& \|\beta_{n}\pi_{n} + (1 - \beta_{n})u_{n} - \pi_{n} +\beta_{n}\pi_{n}- \beta_{n}\pi_{n} \| \\
 \nonumber
   &=& (1 - \beta_{n})\|u_{n} - \pi_{n}\| \\
 \nonumber
&\leq& \|u_{n} - \pi_{n}\| \\
\end{eqnarray}
Again, from  {\color{blue}(\ref{3.7})}, we can see that
$$\lim_{n\to\infty} \|\phi_{n} - \pi_{n}\| = 0$$
\begin{eqnarray}
 \nonumber
  \|\phi_{n} - \delta_{n}\| &=& \| \phi_{n} -  \pi_{n} + \pi_{n} - \delta_{n}\|\\
\nonumber
   &\leq&  \| \phi_{n} -  \pi_{n}\| + \|\pi_{n} - \delta_{n}\|
\end{eqnarray}
Hence $\displaystyle \lim_{n\to\infty} \|\phi_{n} - \delta_{n}\| = 0$
\quad\\
\quad\\
Now from lemma  {\color{blue}(\ref{N})}, lemma  {\color{blue}(\ref{2.1})} and  {\color{blue}(\ref{3})}, we have the following
\begin{eqnarray}\label{3.90}
 \nonumber 
  \|\psi_{n+1} - q\|^{2} &\leq & \|(\alpha_{n}\gamma\varphi(\psi_{n}) + \mu_{n}\xi_{n} + ((1 - \mu_{n})(I - \eta\alpha_{n}\Phi))\psi_{n})- q\|^{2}\\
\nonumber
&=& \|\alpha_{n}(\gamma\varphi(\psi_{n}) - \eta\Phi q) + \mu_{n}(\xi_{n}  - q) + ((1 - \mu_{n})(I - \eta\alpha_{n}\Phi))(\psi_{n}- q)\|^{2}\\
\nonumber
 &\leq & ((1 - \mu_{n}) (1 - \alpha_{n}\tau))^{2}\|\psi_{n} - q\|^{2} + \mu_{n}^{2}\|\xi_{n} - q\|^{2} + \alpha_{n}^{2}\|\gamma\varphi(\psi_{n}) - \eta\Phi q\|^{2}  \\
 \nonumber
&+& 2\mu_{n}\langle ((1 - \mu_{n})( I - \eta\alpha_{n} \Phi))(\psi_{n} - q), \xi_{n} - q\rangle\\
\nonumber
   &+& 2\alpha_{n}\langle \gamma \varphi(\psi_{n}) - \eta\Phi q, \psi_{n+1} - q\rangle \\
 \nonumber
 &\leq & ((1 - \mu_{n}) (1 - \alpha_{n}\tau))^{2}\|\psi_{n} - q\|^{2} + \mu_{n}^{2}\|\xi_{n} - q\|^{2} + \alpha_{n}^{2}\|\gamma\varphi(\psi_{n}) - \eta\Phi q\|^{2}  \\
 \nonumber
&+& 2\mu_{n}((1 - \mu_{n})( 1 - \alpha_{n} \tau))\|\psi_{n} - q\| \|\xi_{n} - q\| + 2\alpha_{n}\|\gamma \varphi(\psi_{n}) - \eta\Phi q\| \|\psi_{n+1} - q\| \\
\nonumber
&\leq&((1 - \mu_{n}) (1 - \alpha_{n}\tau))^{2}\|\psi_{n}- q)\|^{2} + \mu_{n}^{2}\|J_{\lambda_{n}}^{\Pi}(1 - \lambda_{n}\Lambda)\psi_{n}  - J_{\lambda_{n}}^{\Pi}(1 - \lambda_{n}\Lambda)q\|^{2}\\
 \nonumber
 &+&\alpha_{n}^{2}\|\gamma\varphi(\psi_{n}) - \eta\Phi q\|^{2}  \\
  \nonumber
&+&2\mu_{n}((1 - \mu_{n})( 1 - \alpha_{n} \tau)) \|\psi_{n} - q\|^{2} + 2\alpha_{n}\| \gamma\varphi(q) - \eta\Phi q\|\|\psi_{n+1} - q\| \\
\nonumber
&\leq&(1 - 2\alpha_{n}\tau + 2\alpha_{n}\mu_{n}\tau - 2\alpha_{n}^{2}\tau^{2} + \alpha_{n}^{2}\tau^{2}\mu_{n}^{2})\|\psi_{n}- q)\|^{2} \\
\nonumber
&+&  \mu_{n}^{2}[\|\psi_{n}  - q\|^{2} + a(b - 2\alpha)\|\Lambda \psi_{n} - \Lambda q\|^{2}]\\
 \nonumber
 &+&\alpha_{n}^{2}\|\gamma\varphi(\psi_{n}) - \eta\Phi q\|^{2} + 2\alpha_{n}\| \gamma\varphi(q) - \eta\Phi q\|\|\psi_{n+1} - q\| \\
\nonumber
&\leq&(1 - 2\alpha_{n}\tau + 2\alpha_{n}\mu_{n}\tau - 2\alpha_{n}^{2}\tau^{2} + \alpha_{n}^{2}\tau^{2}\mu_{n}^{2} +  \mu_{n}^{2})\|\psi_{n}- q)\|^{2}\\
 \nonumber
 &+& \mu_{n}^{2}a(b - 2\alpha)\|\Lambda \psi_{n} - \Lambda q\|^{2}\\
 \nonumber
 &+&\alpha_{n}^{2}\|\gamma\varphi(\psi_{n}) - \eta\Phi q\|^{2}  + 2\alpha_{n}\| \gamma\varphi(q) - \eta\Phi q\|\|\psi_{n+1} - q\| \\
\nonumber
&\leq&(1 - 2\alpha_{n}\tau)\|\psi_{n}- q\|^{2} -  \mu_{n}^{2}a(2\alpha - b)\|\Lambda \psi_{n} - \Lambda q\|^{2}\\
 \nonumber
 &+&\alpha_{n}^{2}\|\gamma\varphi(\psi_{n}) - \eta\Phi q\|^{2}  + 2\alpha_{n}\| \gamma\varphi(q) - \eta\Phi q\|\|\psi_{n+1} - q\| \\
\nonumber
\end{eqnarray}
Therefore, from  {\color{blue}(\ref{3.90})}, we have
\begin{eqnarray*}
  \mu_{n}^{2} a(2\alpha - b)\|\Lambda \psi_{n} - \Lambda q\|^{2} &\leq & \|\psi_{n}- q\|^{2} - \|\psi_{n+1}- q)\|^{2} -  2\alpha_{n}\tau\|\psi_{n}- q)\|^{2} \\
 &+&\alpha_{n}^{2}\|\gamma\varphi(\psi_{n}) - \eta\Phi q\|^{2}  \\
  &+& 2\alpha_{n}\| \gamma\varphi(q) - \eta\Phi q\|\|\psi_{n+1} - q\| \\
\end{eqnarray*}
Since, $\displaystyle \lim_{n\to\infty}\alpha_{n} = 0$, and from the inequality  {\color{blue}(\ref{3.8})}, with the fact that $\displaystyle \{\psi_{n}\}$ is bounded, we have
\begin{equation}\label{3.10}
    \lim_{\to\infty}\|\Lambda \psi_{n} - \Lambda q\|^{2} = 0
\end{equation}

Since $\displaystyle J_{\lambda_{n}}^{\Pi}$ is $1-$inverse strongly monotone,and for the fact that $\displaystyle  \|\xi_{n} - q\|\leq  \|\delta_{n} - q\|$ then,  this gives us the following
\begin{eqnarray*}
  \|\xi_{n} - q\|^{2} &\leq& \|J_{\lambda_{n}}^{\Pi}(I - \lambda_{n}\Lambda)\psi_{n} - J_{\lambda_{n}}^{\Pi}(I - \lambda_{n}\Lambda )q\|^{2}\\
   &\leq&\langle \xi_{n} - q, (I - \lambda_{n}\Lambda)\psi_{n} - (I - \lambda_{n}\Lambda)q\rangle  \\
   &=& \frac{1}{2}\Big[\|(I - \lambda_{n}\Lambda)\psi_{n} - (I - \lambda_{n}\Lambda) q\|^{2}\\
   &+& \|\xi_{n} - q\|^{2} - \|(I - \lambda_{n}\Lambda)\psi_{n} - (I - \lambda_{n}\Lambda)q - (\xi_{n} - q)\|^{2}\Big] \\
&\leq& \frac{1}{2}\Big[\|\psi_{n} - q\|^{2} + \|\xi_{n} - q\|^{2} - \|\psi_{n} - \xi_{n}\|^{2}\\
   &+& 2\lambda_{n}\langle \xi_{n} - q, \Lambda \psi_{n} - \Lambda q\rangle - \lambda_{n}^{2}\|\Lambda \psi_{n} - \Lambda q \|^{2}\Big] \\
&\leq& \|\psi_{n} - q\|^{2} - \|\psi_{n} - \xi_{n}\|^{2} + 2\lambda_{n}\langle \xi_{n} - q, \Lambda \psi_{n} - \Lambda q\rangle - \lambda_{n}^{2}\|\Lambda \psi_{n} - \Lambda q \|^{2}\\
\end{eqnarray*}
This gives us
\begin{equation}\label{3.12}
  \|\xi_{n} - q\|^{2}\leq \|\psi_{n} - q\|^{2} - \|\psi_{n} - \xi_{n}\|^{2} + 2\lambda_{n}\langle \xi_{n} - q, \Lambda \psi_{n} - \Lambda q\rangle - \lambda_{n}^{2}\|\Lambda \psi_{n} - \Lambda q\|^{2}
\end{equation}
Thus
\begin{eqnarray}\label{3.99}
 \nonumber 
  \|\psi_{n+1} - q\|^{2} &\leq & \|(\alpha_{n}\gamma\varphi(\psi_{n}) + \mu_{n}\xi_{n} + ((1 - \mu_{n})(I - \eta\alpha_{n}\Phi))\psi_{n})- q\|^{2}\\
\nonumber
&=& \|\alpha_{n}(\gamma\varphi(\psi_{n}) - \eta\Phi q) + \mu_{n}(\xi_{n}  - q) + ((1 - \mu_{n})(I - \eta\alpha_{n}\Phi)(\psi_{n}- q))\|^{2}\\
\nonumber
 &\leq & ((1 - \mu_{n}) (1 - \alpha_{n}\tau))^{2}\|\psi_{n} - q\|^{2} + \mu_{n}^{2}\|\xi_{n} - q\|^{2} + \alpha_{n}^{2}\|\gamma\varphi(\psi_{n}) - \eta\Phi q\|^{2}  \\
 \nonumber
&+& 2\mu_{n}\langle ((1 - \mu_{n})( I - \eta\alpha_{n} \Phi))(\psi_{n} - q), \xi_{n} - q\rangle\\
\nonumber
   &+& 2\alpha_{n}\langle \gamma \varphi(\psi_{n}) - \eta\Phi q, \psi_{n+1} - q\rangle \\
 \nonumber
 &\leq & ((1 - \mu_{n}) (1 - \alpha_{n}\tau))^{2}\|\psi_{n} - q\|^{2} + \mu_{n}^{2}\|\xi_{n} - q\|^{2} + \alpha_{n}^{2}\|\gamma\varphi(\psi_{n}) - \eta\Phi q\|^{2}  \\
 \nonumber
&+& 2\mu_{n}((1 - \mu_{n})( 1 - \alpha_{n} \tau))\|\psi_{n} - q\| \|\xi_{n} - q\| + 2\alpha_{n}\|\gamma \varphi(\psi_{n}) - \eta\Phi q\| \|\psi_{n+1} - q\| \\
 \nonumber
&\leq&((1 - \mu_{n}) (1 - \alpha_{n}\tau))^{2}\|\psi_{n}- q)\|^{2} + \mu_{n}^{2}\Big[\|\psi_{n} - q\|^{2} - \|\psi_{n} - \xi_{n}\|^{2}\\
\nonumber
 &+& 2\lambda_{n}\langle \xi_{n} - q, \Lambda \psi_{n} - \Lambda q\rangle - \lambda_{n}^{2}\|\Lambda \psi_{n} - \Lambda q\|^{2}\Big] + \alpha_{n}^{2}\|\gamma\varphi(\psi_{n}) - \eta\Phi q\|^{2}\\
 \nonumber
&+& 2\mu_{n}((1 - \mu_{n})( 1 - \alpha_{n} \tau))\|\psi_{n} - q\| \|\xi_{n} - q\| + 2\alpha_{n}\|\gamma \varphi(\psi_{n}) - \eta\Phi q\| \|\psi_{n+1} - q\| \\
\nonumber
&\leq&(1 - 2\alpha_{n}\tau)\|\psi_{n}- q)\|^{2} - \mu_{n}^{2}\|\psi_{n} - \xi_{n}\|^{2}\\
\nonumber
 &+& 2\mu_{n}^{2}\lambda_{n}\langle \xi_{n} - q, \Lambda \psi_{n} - \Lambda q\rangle - \lambda_{n}^{2}\|\Lambda \psi_{n} - \Lambda q\|^{2}\\
  &+& \alpha_{n}^{2}\|\gamma\varphi(\psi_{n}) - \eta\Phi q\|^{2} + 2\alpha_{n}\|\gamma \varphi(\psi_{n}) - \eta\Phi q\| \|\psi_{n+1} - q\| \\
\nonumber
\end{eqnarray}
Therefor, since $\displaystyle \alpha_{n}\to 0$ as $\displaystyle n\to \infty$ with inequalities {\color{blue}\ref{3.8}} and {\color{blue}}\ref{3.10}, we have
$$\lim_{n\to\infty}\|\psi_{n} - \xi_{n}\| = 0$$

From  {\color{blue}(\ref{3.0})} and lemma  {\color{blue}(\ref{2.2})} with the fact that $\displaystyle T_{3}$ is quasi-nonexpansive, we have the following estimate
\begin{eqnarray}\label{3.11}
 \nonumber
  \|\xi_{n} - q\|^{2} &=& \|\gamma_{n}\phi_{n} + (1 - \gamma_{n})z_{n} - q\|^{2}  \\
\nonumber
   &=&\gamma_{n}\|\phi_{n} - q\|^{2} + (1 - \gamma_{n})\|z_{n} - q\|^{2} -   (1 - \gamma_{n})\gamma_{n}\|z_{n} - \phi_{n}\|^{2} \\
\nonumber
 &=&\gamma_{n}\|\phi_{n} - q\|^{2} + (1 - \gamma_{n})H(T_{3}\phi_{n}, T_{3}q)^{2} -   (1 - \gamma_{n})\gamma_{n}\|z_{n} - \phi_{n}\|^{2} \\
\nonumber
 &=&\gamma_{n}\|\phi_{n} - q\|^{2} + (1 - \gamma_{n})\|\phi_{n} - q \|^{2} -   (1 - \gamma_{n})\gamma_{n}\|z_{n} - \phi_{n}\|^{2} \\
  &\leq& \|\psi_{n} - q \|^{2} -  (1 - \gamma_{n})\gamma_{n}\|z_{n} - \phi_{n}\|^{2}
\end{eqnarray}
Therefore
\begin{eqnarray}\label{3.9}
 \nonumber 
  \|\psi_{n+1} - q\|^{2} &\leq & \|(\alpha_{n}\gamma\varphi(\psi_{n}) + \mu_{n}\xi_{n} + ((1 - \mu_{n})(I - \eta\alpha_{n}\Phi))\psi_{n})- q\|^{2}\\
\nonumber
&=& \|\alpha_{n}(\gamma\varphi(\psi_{n}) - \eta\Phi q) + \mu_{n}(\xi_{n}  - q) + ((1 - \mu_{n})(I - \eta\alpha_{n}\Phi)(\psi_{n}- q))\|^{2}\\
\nonumber
 &\leq & ((1 - \mu_{n}) (1 - \alpha_{n}\tau))^{2}\|\psi_{n} - q\|^{2} + \mu_{n}^{2}\|\xi_{n} - q\|^{2} + \alpha_{n}^{2}\|\gamma\varphi(\psi_{n}) - \eta\Phi q\|^{2}  \\
 \nonumber
&+& 2\mu_{n}\langle ((1 - \mu_{n})( I - \eta\alpha_{n} \Phi))(\psi_{n} - q), \xi_{n} - q\rangle\\
\nonumber
   &+& 2\alpha_{n}\langle \gamma \varphi(\psi_{n}) - \eta\Phi q, \psi_{n+1} - q\rangle \\
 \nonumber
 &\leq & ((1 - \mu_{n}) (1 - \alpha_{n}\tau))^{2}\|\psi_{n} - q\|^{2} + \mu_{n}^{2}\|\xi_{n} - q\|^{2} + \alpha_{n}^{2}\|\gamma\varphi(\psi_{n}) - \eta\Phi q\|^{2}  \\
 \nonumber
&+& 2\mu_{n}((1 - \mu_{n})( 1 - \alpha_{n} \tau))\|\psi_{n} - q\| \|\xi_{n} - q\| + 2\alpha_{n}\|\gamma \varphi(\psi_{n}) - \eta\Phi q\| \|\psi_{n+1} - q\| \\
 \nonumber
 &\leq & ((1 - \mu_{n}) (1 - \alpha_{n}\tau))^{2}\|\psi_{n} - q\|^{2} + \mu_{n}^{2}\|\psi_{n} - q \|^{2} -  \mu_{n}^{2}(1 - \gamma_{n})\gamma_{n}\|z_{n} - \phi_{n}\|^{2}\\
 \nonumber
 & +& \alpha_{n}^{2}\|\gamma\varphi(\psi_{n}) - \eta\Phi q\|^{2}  +  2\mu_{n}((1 - \mu_{n})( 1 - \alpha_{n} \tau))\|\psi_{n} - q\| \|\xi_{n} - q\|\\
  \nonumber
  &+& 2\alpha_{n}\|\gamma \varphi(\psi_{n}) - \eta\Phi q\| \|\psi_{n+1} - q\| \\
 \nonumber
 &\leq & (1 - 2\alpha_{n}\tau)\|\psi_{n} - q\|^{2} -  \mu_{n}^{2}(1 - \gamma_{n})\gamma_{n}\|z_{n} - \phi_{n}\|^{2}\\
 \nonumber
 & +& \alpha_{n}^{2}\|\gamma\varphi(\psi_{n}) - \eta\Phi q\|^{2}  + 2\alpha_{n}\|\gamma \varphi(\psi_{n}) - \eta\Phi q\| \|\psi_{n+1} - q\| \\
 \nonumber
&\leq & \|\psi_{n} - q\|^{2} - 2\alpha_{n}\tau\|\psi_{n} - q\|^{2} -  \mu_{n}^{2}(1 - \gamma_{n})\gamma_{n}\|z_{n} - \phi_{n}\|^{2}\\
 \nonumber
 & +& \alpha_{n}^{2}\|\gamma\varphi(\psi_{n}) - \eta\Phi q\|^{2}  + 2\alpha_{n}\|\gamma \varphi(\psi_{n}) - \eta\Phi q\| \|\psi_{n+1} - q\| \\
 \nonumber
\end{eqnarray}

Now setting $\displaystyle \mathbb{M} =  \alpha_{n}\|\gamma\varphi(\psi_{n}) - \eta\Phi q\|^{2} - 2\tau\|\psi_{n} - q\|^{2} + 2\|\gamma \varphi(\psi_{n}) - \eta\Phi q\| \|\psi_{n+1} - q\| > 0$,and for the fact that the sequence $\displaystyle \{\psi_{n}\}$ is bounded such that
\begin{equation}\label{3.15}
  \mu_{n}^{2}(1 - \gamma_{n})\gamma_{n}\|z_{n} - \phi_{n}\|^{2} \leq  \|\psi_{n}- q\|^{2} - \|\psi_{n + 1}- q)\|^{2}  + \alpha_{n}\mathbb{M}
\end{equation}
Therefore, from  {\color{blue}(\ref{3.15})} and  {\color{blue}(\ref{3.8})}, we have
\begin{equation}\label{3.16}
    \lim_{n\to\infty}(1 - \gamma_{n})\gamma_{n}\|z_{n} - \phi_{n}\|^{2} = 0
\end{equation}
Since $\displaystyle \lim_{n\to\infty}\inf((1 - \gamma_{n})\gamma_{n}) > 0$
\begin{equation}\label{3.17}
    \lim_{n\to\infty}\|z_{n} - \phi_{n}\| = 0
\end{equation}
Again, with $\displaystyle z \in T_{3}\phi_{n}$
\begin{equation}\label{3.18}
    \lim_{n\to\infty}d(\phi_{n},T_{3}\phi_{n}) = 0
\end{equation}
Moreover, since $H$ is reflexive and $\displaystyle \{\psi_{n}\}$ is bounded, we then prove that $\displaystyle \lim_{n\to +\infty}\sup\langle \eta \Phi\psi^{*} - \gamma\varphi(\psi^{*}), \psi^{*} - \psi_{n}\rangle \leq 0$. We let the subsequence $\displaystyle \{\psi_{n_i}\}$ of   $\displaystyle \{\psi_{n}\}$  to converge weakly to  $\displaystyle \psi^{**}$ in $\mathbb{K}$, and
$$\displaystyle \lim_{n\to +\infty}\langle \eta \Phi\psi^{*} - \gamma\varphi(\psi^{*}), \psi^{*} - \psi_{n}\rangle = \lim_{n\to +\infty}\langle \eta \Phi\psi^{*} - \gamma\varphi(\psi^{*}), \psi^{*} - \psi_{n_ i}\rangle$$
Again, since $\displaystyle I - T_{1},  I - T_{2}$ and  $\displaystyle I - T_{3}$ satisfies the demiclosed principle and from  {\color{blue}(\ref{333})},  {\color{blue}(\ref{334})} and  {\color{blue}(\ref{3.18})}, we obtain $\displaystyle \psi^{**}\in F ix(T_{1})\cap F ix(T_{2})\cap Fix(T_{3})$. We now show that $\displaystyle \psi^{**}\in S(\Pi, \Lambda)$. Since $\Lambda$ is $\alpha-$inverse strongly monotone, $\Lambda$ is Lipschitz continuous mapping. From lemma  {\color{blue}(\ref{5.0})}, it follows that $(\Pi + \Lambda)$ is maximal monotone.
\quad\\
Let $(\nu,g)\in G(\Pi + \Lambda)$, that is $g - \Lambda \nu \in \Pi(\nu)$. Since $\displaystyle \delta_{n_i} = J_{\lambda_{n_i}}^{\Pi}(\psi_{n_i}) - \lambda_{n_i}\Lambda\psi_{n_i})$, we have $\displaystyle \psi_{n_i} - \lambda_{n_i}\psi_{n_i}\in (I + \lambda_{n_i}\Pi)\delta_{n_i}$, that is
$\displaystyle \frac{1}{\lambda_{n_i}}(\psi_{n_i} - \delta_{n_i} - \lambda_{n_i}\Lambda\psi_{n_i})\in\Pi(\delta_{n_i}).$ By maximal monotonocity of $(\Pi + \Lambda)$, gives $$\langle \nu - \delta_{n_i}, g - \Lambda \nu -  \frac{1}{\lambda_{n_i}}(\psi_{n_i} - \delta_{n_i} - \lambda_{n_i}\Lambda\psi_{n_i})\geq 0$$ and, therefore
\begin{eqnarray}
 \nonumber
  \langle \nu - \delta_{n_i}, g\rangle&\geq & \langle \nu - \delta_{n_i},\Lambda \nu - \frac{1}{\lambda_{n_i}}(\psi_{n_i} - \delta_{n_i} - \lambda_{n_i}\Lambda\psi_{n_i})\rangle \\
   \nonumber
   &=&\langle \nu - \delta_{n_i}, \Lambda \nu - \Lambda\delta_{n_i} + \Lambda\delta_{n_i} +  \frac{1}{\lambda_{n_i}}(\psi_{n_i} - \delta_{n_i} - \lambda_{n_i}\Lambda\psi_{n_i})\rangle \\
    \nonumber
   &\geq& \langle \nu - \delta_{n_i}, \Lambda \nu - \Lambda\psi_{n_i}\rangle + \langle \nu - \delta_{n_i},\frac{1}{\lambda_{n_i}}(\psi_{n_i} - \delta_{n_i}\rangle
\end{eqnarray}
It then follows from $\displaystyle \|\delta_{n} - \psi_{n}\| \to 0, \|\Lambda\delta_{n} - \Lambda\psi_{n}\| \to 0\quad and \quad \delta_{n_i} \to \psi^{**}$ weakly that $\displaystyle \lim_{n\to\infty }\langle \nu - \delta_{n_i}, g\rangle = \langle \nu - \psi^{**}, g\rangle$ and hence $\displaystyle \psi^{**} \in S(\Pi, \Lambda)$. Therefore, $\displaystyle \psi^{**}\in\Omega$,

\quad\\
On the other hand, for the fact that $\psi^{*}$ solves the variational inequality  {\color{blue}(\ref{331})}.

\begin{eqnarray}
 \nonumber
  \lim_{n\to +\infty}\sup\langle \eta\Phi\psi^{*} - \gamma\varphi(\psi^{*}), \psi^{*} - \psi_{n}\rangle &=&\lim_{n\to +\infty}\sup\langle \eta\Phi\psi^{*} - \gamma\varphi(\psi^{*}), \psi^{*} - \psi_{n_i}\rangle  \\
   \nonumber
   &=& \langle \eta\Phi\psi^{*} - \gamma\varphi(\psi^{*}), \psi^{*} - \psi^{**}\rangle\leq 0 \\
\end{eqnarray}
Lastly, we now prove that $\displaystyle \lim_{n\to\infty}\|\psi_{n} - \psi^{*}\| = 0$,  that is $\displaystyle \psi_{n} \to  \psi^{*} \quad as \quad n \to \infty$. From lemma  {\color{blue}(\ref{2.1})},  {\color{blue}(\ref{3})} and  {\color{blue}(\ref{3.11})}, we have the following
\begin{eqnarray}\label{99}
 \nonumber
  \|\psi_{n+1} - \psi^{*}\|^{2} &=&  \|P_{\mathbb{K}}(\alpha_{n}\gamma\varphi(\psi_{n}) + \mu_{n}\xi_{n} + ((1 - \mu_{n})(I - \eta\alpha_{n}\Phi))\psi_{n})- \psi^{*}\|^{2}\\
\nonumber
&\leq&  \langle\alpha_{n}\gamma\varphi(\psi_{n}) + \mu_{n}\xi_{n} + ((1 - \mu_{n})(I - \eta\alpha_{n}\Phi))\psi_{n}- \psi^{*}, \psi_{n+1} - \psi^{*}\rangle\\
\nonumber
&\leq&  \langle\alpha_{n}\gamma\varphi(\psi_{n}) - \alpha_{n}\gamma\varphi(\psi^{*}) + \mu_{n}\xi_{n}  + \mu_{n}\xi^{*} + ((1 - \mu_{n})(I - \eta\alpha_{n}\Phi))\psi_{n}- \psi^{*} \\
\nonumber
&+& \alpha_{n}\gamma\varphi(\psi^{*}) - \mu_{n}\xi^{*} - \alpha_{n}\eta\Phi\psi^{*} + \alpha_{n}\eta\Phi\psi^{*} , \psi_{n+1} - \psi^{*}\rangle\\
\nonumber
&\leq& \alpha_{n}\gamma \|\varphi(\psi_{n})  - \varphi(\psi^{*})\|\|\psi_{n+1} - \psi^{*}\| + \mu_{n}\|\xi_{n}  - \psi^{*}\|\|\psi_{n+1} - \psi^{*}\| \\
\nonumber
&+& \|((1 - \mu_{n})(I - \eta\alpha_{n}\Phi))(\psi_{n}- \psi^{*})\|\|\psi_{n+1} - \psi^{*}\| + \alpha_{n}\langle \eta\Phi\psi^{*} - \gamma\varphi(\psi^{*}), \psi_{n+1} - \psi^{*}\rangle\\
\nonumber
 &\leq& \alpha_{n}\gamma b\|\psi_{n}  - \psi^{*}\|\|\psi_{n+1} - \psi^{*}\| + \mu_{n}\|\psi_{n}  - \psi^{*}\|\|\psi_{n+1} - \psi^{*}\| \\
\nonumber
&+& ((1 - \mu_{n})(1 - \tau\alpha_{n}))\|\psi_{n}- \psi^{*}\|\|\psi_{n+1} - \psi^{*}\| + \alpha_{n}\langle \eta\Phi\psi^{*} - \gamma\varphi(\psi^{*}), \psi_{n+1} - \psi^{*}\rangle\\
\nonumber
 &\leq& [1 - \alpha_{n}( \tau - \gamma b - \mu_{n}\tau)]\|\psi_{n}  - \psi^{*}\|\|\psi_{n+1} - \psi^{*}\|  + \alpha_{n}\langle \eta\Phi\psi^{*} - \gamma\varphi(\psi^{*}), \psi_{n+1} - \psi^{*}\rangle\\
\nonumber
&\leq& [1 - \alpha_{n}( \tau - \gamma b - \mu_{n}\tau)]\|\psi_{n}  - \psi^{*}\|^{2} + 2\alpha_{n}\langle \eta\Phi\psi^{*} - \gamma\varphi(\psi^{*}), \psi_{n+1} - \psi^{*}\rangle\\
\nonumber
 \nonumber
\end{eqnarray}
Thus, from lemma  {\color{blue}(\ref{29})}, it follows that $\displaystyle \psi_{n} \to \psi^{*} \quad\quad as\quad\quad n \to \infty.$
  \item[Case 2]
  Suppose that the sequence $\displaystyle \Big\{\|\psi_{n} - \psi^{*}\|\Big\}$ is monotonically increasing. Set $\mathbb{B}_{n}: = \|\psi_{n} - \psi^{*}\|^{2}$ and $\tau: = \mathbb{N}\to\mathbb{N}$ be a mapping for all $n\geq n_{0}$ (for some $n_{0}$ sufficient large), by $\displaystyle \tau_{n}: = \max\{k\in\mathbb{N}: k\leq n, \mathbb{B}_{k}\leq\mathbb{B}_{k+1}\}$. Then, $\tau$ is a nondecreasing sequence, such that $\tau_{n} \to \infty$ as $n\to\infty$ and $\displaystyle \mathbb{B}_{\tau(n)}\leq\mathbb{B}_{\tau(n)+1}\}$ for all $n\geq n_{0}$. Now, from  {\color{blue}(\ref{3.50})}, we have
  \begin{eqnarray}\label{3.5}
\nonumber
&&\mu_{\tau (n)}(1 - \alpha_{\tau (n)}\tau(1 - \mu_{\tau (n)}))\Big[ (1 - \theta_{\tau (n)})(\theta_{\tau (n)} - \beta)\|v_{\tau (n)} - \delta_{\tau (n)}  \|^{2}\\
\nonumber
 &+&(1 - \beta_{\tau (n)})(\beta_{\tau (n)} - \beta)\|u_{\tau (n)} - \pi_{\tau (n)}  \|^{2}\Big]\\
\nonumber
   &\leq& \|\psi_{\tau (n)} - q\|^{2} - \|\psi_{\tau (n)+ 1} - q\|^{2}\\
\nonumber
&-& \alpha_{n}\tau(\mu_{\tau (n)} - \mu_{\tau (n)}^{2}   -  \alpha_{\tau (n)}\tau)\|\psi_{\tau (n)} - q\|^{2} +  2\alpha_{\tau (n)}M .\\
\end{eqnarray}
$$\mu_{\tau (n)}(1 - \alpha_{\tau (n)}\tau(1 - \mu_{\tau (n)}))\Big[ (1 - \theta_{\tau (n)})(\theta_{\tau (n)} - \beta)\|v_{\tau (n)} - \delta_{\tau (n)}  \|^{2}+ (1 - \beta_{\tau (n)})(\beta_{\tau (n)} - \beta)\|u_{\tau (n)} - \pi_{\tau (n)}  \|^{2}\Big] = 0$$
  Since $\displaystyle \beta_{\tau(n)}\in (\beta,1) \quad and \quad \theta_{\tau(n)}\in (\beta,1)$ and $$\displaystyle \lim_{n\to \infty}\inf\beta_{\tau(n)}(1 - \beta_{\tau(n)})> 0\quad and \quad \lim_{n\to \infty}\inf\theta_{\tau(n)}(1 - \theta_{\tau(n)})> 0,$$ we have
  \begin{equation*}
    \lim_{n\to\infty}\|u_{\tau(n)} - \pi_{\tau(n)}\| = 0 \quad\quad and \quad\quad  \lim_{n\to\infty}\|v_{\tau(n)} - \delta_{\tau(n)}\| = 0
  \end{equation*}
 With $\displaystyle v_{\tau(n)}\in T_{1} \delta_{\tau(n)}$ and $\displaystyle u_{\tau(n)}\in T_{2} \pi_{\tau(n)}$, it follows that
  \begin{equation*}
    \lim_{n\to\infty}d\Big(\delta_{\tau(n)}, T_{1} \delta_{\tau(n)}\Big) = 0 \quad\quad and \quad\quad   \lim_{n\to\infty}d\Big(\pi_{\tau(n)}, T_{2} \pi_{\tau(n)}\Big) = 0
  \end{equation*}

  Following the same argument in case $1$, we conclude that
  $$\lim_{\tau(n)\to + \infty}\sup\langle \eta\Phi\psi^{*} - \gamma\varphi(\psi^{*}), \psi^{*} - \psi_{\tau(n) +1}\rangle \leq 0$$
  Therefore, for all $\displaystyle n\geq n_{0}$, we have
  \begin{eqnarray*}
  0\leq \|\psi_{\tau(n)+1} - \psi^{*}\|^{2}-  \|\psi_{\tau(n)} - \psi^{*}\|^{2} &\leq&\alpha_{\tau(n)}\Big[ -( \tau - \gamma b - \mu_{n}\tau)\|\psi_{n}  - \psi^{*}\|^{2} \\
  &+& 2\langle \eta\Phi\psi^{*} - \gamma\varphi(\psi^{*}), \psi_{n+1} - \psi^{*}\rangle \Big]\\
\end{eqnarray*}
  Hence, we have
  \begin{equation*}
  \|\psi_{\tau(n)} - \psi^{*}\|^{2} \leq\frac{2}{\tau - \gamma b - \mu_{n}\tau}\langle \gamma\varphi(\psi_{n}) - \eta\Phi\psi^{*}, \psi_{\tau(n)+1} - \psi^{*}\rangle
\end{equation*}
Then we have $\displaystyle \lim_{n\to\infty}\|\psi_{\tau(n)} - \psi^{*}\|^{2} = 0$. Therefore $\displaystyle \lim_{n\to \infty}\mathbb{B}_{\tau(n)} = \lim_{n\to \infty}\mathbb{B}_{\tau(n) + 1} = 0$
\quad\\
Furthermore, for all $n\geq n_{0}$, we have $\displaystyle \mathbb{B}_{\tau(n)} \leq \mathbb{B}_{\tau(n) + 1}$ if $\displaystyle n \neq \tau(n)$ (that is $\displaystyle n >\tau(n)$, because $\displaystyle \mathbb{B}_{j} > \mathbb{B}_{j + 1}, \quad\quad for\quad\quad \tau(n+1) \leq j \leq n$
\quad\\
Hence, $\displaystyle 0 \leq \mathbb{B}_{\tau(n)} \leq \max\Big\{\mathbb{B}_{\tau(n)}, \mathbb{B}_{\tau(n) + 1}\Big\} = \mathbb{B}_{\tau(n) + 1}$. Therefore, $\displaystyle \mathbb{B}_{n} \to 0, \quad as \quad n \to\infty$ and this implies that $\displaystyle \psi_{n}\to \psi^{*}\quad as \quad n\to\infty.$ This complete the proof.
\end{description}
\end{proof}
Now using the same argument of the proof of theorem  {\color{blue}(\ref{NB1})} in theorem  {\color{blue}(\ref{NB2})}, and theorem  {\color{blue}(\ref{NB3})}, we achieved the desired results. In theorem  {\color{blue}(\ref{NB2})}, the demiclosed principle assumption is not considered and finally in theorem  {\color{blue}(\ref{NB3})}, we let $\displaystyle T_{1} = P_{T_1},  T_{2} = P_{T_2}$ and  $\displaystyle T_{3} = P_{T_3}$ without the assumptions that $T_{1}q = T_{2}q = T_{3}q = \{q\}, \forall q \in \Omega$
\begin{theorem}\label{NB2}
Let $H$ be a real Hilbert space and $\mathbb{K}$ be a nonempty, closed convex subset of $H$. Let $\Lambda : \mathbb{K} \to H$ be an $\alpha-$inverse strongly monotone operator and let $\Phi : H \to H$ be an $k-$strongly monotone and $L-$Lipschitzian operator. Let $\varphi : \mathbb{K} \to H$ be an $b-$Lipschitzian mapping and $\Pi : H \to 2^{H}$ be a maximal monotone mapping such that the domain of $\Pi$ is included in $\mathbb{K}$. Let $T_{1},T_{2} : \mathbb{K}\to CB(\mathbb{K})$ be a multivalued $\beta-$ demicontractive mapping and $T_{3} : \mathbb{K} \to CB(\mathbb{K})$ be a multivalued quasi-nonexpansive mapping such that $\Omega := F ix(T_{1})\cap F ix(T_{2})\cap Fix(T_{3})\cap S(\Pi, \Lambda) \neq \emptyset$
and $T_{1}q = T_{2}q = T_{3}q = \{q\}, \forall q \in \Omega$. For given $\psi_{0} \in \mathbb{K}$, let $\{\psi_{n}\}$ be generated by the algorithm:

\begin{equation}\label{33}
    \left\{
      \begin{array}{ll}
        \delta_{n} = J_{\lambda_{n}}^{\Pi}(I - \lambda_{n}\Lambda)\psi_{n};\\
        \\
        \pi_{n} = \theta_{n}\delta_{n} + (1 - \theta_{n})v_{n}, \quad v_{n} \in T_{1}\delta_{n}; \\
        \\
        \phi_{n} = \beta_{n}\pi_{n} + (1 - \beta_{n})u_{n}, \quad u_{n}\in T_{2}\pi_{n}; \\
        \\
        \xi_{n} = \gamma_{n}\phi_{n} + (1 - \gamma_{n})z_{n}, \quad z_{n}\in T_{3}\phi_{n}; \\
        \\
        \psi_{n+1} = P_{\mathbb{K}}(\alpha_{n}\gamma \varphi(\psi_{n}) + \mu_{n}\xi_{n} + ((1 - \mu_{n}) (I - \eta\alpha_{n} \Phi))\psi_{n})
      \end{array}
    \right.
\end{equation}
where $\displaystyle\{\beta_{n}\} , \{\gamma_{n}\},  \{\theta_{n}\},  \{\lambda_{n}\}$ and $\displaystyle \{\alpha_{n}\}$ are real sequence in $(0,1)$ satisfying the following conditions
\begin{description}
  \item[i)] $\displaystyle \lim_{n\to\infty}\alpha_{n} = 0\quad\quad \sum_{n=0}^{\infty}\alpha_{n} < \infty$
  \item[ii)] $\displaystyle \lim_{n\to\infty}\inf(1 - \beta_{n})(\beta_{n} - \beta)> 0\quad and \quad \lim_{n\to\infty}\inf(1 - \theta_{n})(\theta_{n} - \beta)> 0, (\beta_{n},\theta_{n})\in (\beta,1)$
  \item[iii)] $\lim_{n\to\infty}\inf(1 - \gamma_{n})\gamma_{n})> 0 ,\lim_{n\to\infty}\inf(1 - \beta_{n})\beta_{n})> 0 \quad\quad and \quad\quad \lim_{n\to\infty}\inf(1 - \theta_{n})\theta_{n})> 0 $
\end{description}
Assume that $\displaystyle 0 < \eta < \frac{2k}{L^{2}}, 0 < \gamma b < \tau$, where $\displaystyle \tau = \eta\Big(k - \frac{L^{2}\eta}{2}\Big)$, and the sequences defined in {\color{blue}(\ref{33})}, that is $\displaystyle \{\psi_{n}\}$ and $\displaystyle \{\delta_{n}\}$ converge strongly to unique solution $\displaystyle \psi^{*} \in \Omega$, which also solve the following variational inequality:
\begin{equation}\label{34}
    \langle \eta \Phi \psi^{*} - \gamma \phi(\psi^{*}), \psi^{*} - q\rangle \leq 0, \quad \forall q \in \Omega
\end{equation}
\end{theorem}
\quad\\

\begin{theorem}\label{NB3}
Let $H$ be a real Hilbert space and $\mathbb{K}$ be a nonempty, closed convex subset of $H$. Let $\Lambda : \mathbb{K} \to H$ be an $\alpha-$inverse strongly monotone operator and let $\Phi : H \to H$ be an $k-$strongly monotone and $L-$Lipschitzian operator. Let $\varphi : \mathbb{K} \to H$ be an $b-$Lipschitzian mapping and $\Pi : H \to 2^{H}$ be a maximal monotone mapping such that the domain of $\Pi$ is included in $\mathbb{K}$. Let $T_{1},T_{2} : \mathbb{K} \to CB(\mathbb{K})$ be a multivalued $\beta-$ demicontractive mapping and $T_{3} : \mathbb{K} \to CB(\mathbb{K})$ be a multivalued quasi-nonexpansive mapping such that $\Omega := F ix(T_{1})\cap F ix(T_{2})\cap Fix(T_{3})\cap S(\Pi, \Lambda) \neq \emptyset.$
 For given $\psi_{0} \in \mathbb{K}$, let $\{\psi_{n}\}$ be generated by the algorithm:

\begin{equation}\label{336}
    \left\{
      \begin{array}{ll}
        \delta_{n} = J_{\lambda_{n}}^{\Pi}(I - \lambda_{n}\Lambda)\psi_{n};\\
        \\
        \pi_{n} = \theta_{n}\delta_{n} + (1 - \theta_{n})v_{n}, \quad v_{n} \in T_{1}\delta_{n}; \\
        \\
        \phi_{n} = \beta_{n}\pi_{n} + (1 - \beta_{n})u_{n}, \quad u_{n}\in T_{2}\pi_{n}; \\
        \\
        \xi_{n} = \gamma_{n}\phi_{n} + (1 - \gamma_{n})z_{n}, \quad z_{n}\in T_{3}\phi_{n}; \\
        \\
        \psi_{n+1} = P_{\mathbb{K}}(\alpha_{n}\gamma \varphi(\psi_{n}) + \mu_{n}\xi_{n} + ((1 - \mu_{n}) (I - \eta\alpha_{n} \Phi))\psi_{n})
      \end{array}
    \right.
\end{equation}
where $\displaystyle\{\beta_{n}\}, \{\gamma_{n}\},  \{\theta_{n}\},  \{\lambda_{n}\}$ and $\displaystyle \{\alpha_{n}\}$ are real sequence in $(0,1)$ satisfying the following conditions
\begin{description}
  \item[i)] $\displaystyle \lim_{n\to\infty}\alpha_{n} = 0\quad\quad \sum_{n=0}^{\infty}\alpha_{n} < \infty$
  \item[ii)] $\displaystyle \lim_{n\to\infty}\inf(1 - \beta_{n})(\beta_{n} - \beta)> 0\quad and \quad \lim_{n\to\infty}\inf(1 - \theta_{n})(\theta_{n} - \beta)> 0, (\beta_{n},\theta_{n})\in (\beta,1)$
 \item[iii)] $\lim_{n\to\infty}\inf(1 - \gamma_{n})\gamma_{n})> 0 ,\lim_{n\to\infty}\inf(1 - \beta_{n})\beta_{n})> 0 \quad\quad and \quad\quad \lim_{n\to\infty}\inf(1 - \theta_{n})\theta_{n})> 0 $
\end{description}
Assume that $\displaystyle 0 < \eta < \frac{2k}{L^{2}}, 0 < \gamma b < \tau$, where $\displaystyle \tau = \eta\Big(k - \frac{L^{2}\eta}{2}\Big)$, and $I - P_{T_1}, I - P_{T_2}\quad and \quad I - P_{T_3}$ are demiclosed at origin. Hence, the sequences defined in  {\color{blue}(\ref{333})}, that is $\displaystyle \{\psi_{n}\}$ and $\displaystyle \{\delta_{n}\}$ converge strongly to unique solution $\displaystyle \psi^{*} \in \Omega$, which also solve the following variational inequality:
\begin{equation}\label{35}
    \langle \eta \Phi \psi^{*} - \gamma \phi(\psi^{*}), \psi^{*} - q\rangle \leq 0, \quad \forall q \in \Omega
\end{equation}
\end{theorem}
\quad\\

\section{Conclusion}
The modified general viscosity iterative process presented in this research offers a powerful tool for solving variational inclusion and fixed point problems involving and
and fixed point problem with respectively set-valued maximal monotone mapping and inverse strongly monotone and multivalued quasi-nonexpansive and demicontractive operators. By incorporating the method into optimization algorithms, it enhances the optimization process, improves convergence properties, and enables the handling of complex constraints and objective functions. Its application finds relevance in diverse domains, yielding more efficient and accurate optimization solutions.
\section*{Declarations}

The authors hereby declare that there is no potential conflict of interest with respect to this research,or with publication of this article. Further, authors also declare that they have no known competing financial interests or personal relationships that could have appeared to influence the work reported in this article.

\end{document}